\documentclass{IEEEtran}

\usepackage{graphicx}
\usepackage{graphics}
\usepackage{amssymb}
\usepackage{amsmath}
\usepackage{subfigure}
\usepackage{epsfig}
\usepackage{color}
\usepackage{enumitem}
\usepackage{float}
\usepackage{soul}
\usepackage{wrapfig}
\usepackage[noadjust]{cite}
\usepackage{url}
\usepackage{stmaryrd}

\usepackage{stackengine}
\newcommand\barbelow[1]{\stackunder[1.2pt]{$#1$}{\rule{.8ex}{.075ex}}}

\newtheorem{Remark 1}{Remark}
\newtheorem{Remark 2}[Remark 1]{Remark}
\newtheorem{Remark 3}[Remark 1]{Remark}
\newtheorem{Remark 4}[Remark 1]{Remark}
\newtheorem{Remark 5}[Remark 1]{Remark}
\newtheorem{Remark 6}[Remark 1]{Remark}
\newtheorem{Remark 7}[Remark 1]{Remark}

\newtheorem{Theorem 1}{Theorem}
\newtheorem{Theorem 2}[Theorem 1]{Theorem}
\newtheorem{Theorem 3}[Theorem 1]{Theorem}
\newtheorem{Theorem 4}[Theorem 1]{Theorem}
\newtheorem{Theorem 5}[Theorem 1]{Theorem}
\newtheorem{Theorem 6}[Theorem 1]{Theorem}
\newtheorem{Theorem 7}[Theorem 1]{Theorem}
\newtheorem{Theorem 8}[Theorem 1]{Theorem}
\newtheorem{Theorem 9}[Theorem 1]{Theorem}
\newtheorem{Theorem 10}[Theorem 1]{Theorem}

\newtheorem{Lemma 1}{Lemma}
\newtheorem{Lemma 2}[Lemma 1]{Lemma}

\newtheorem{Definition 1}{Definition}

\title{\LARGE \bf Secure and Privacy-Preserving Consensus}

\author{Minghao Ruan, Huan Gao,~\IEEEmembership{Student~Member,~IEEE,} and Yongqiang Wang,~\IEEEmembership{Senior~Member,~IEEE}
\thanks{Part of the results was accepted for presentation at the Cyber-Physical Systems Security and Privacy Workshop \cite{CPS_SPC}. The work was supported in part by the National Science Foundation under Grant 1824014 and 1738902.}
\thanks{Minghao Ruan, Huan Gao, and Yongqiang Wang are with the Department of Electrical and Computer Engineering, Clemson University, Clemson, SC 29634, USA. Corresponding author: Yongqiang Wang ({\tt\small{yongqiw@clemson.edu}}).}}

\begin{document}

\maketitle
\thispagestyle{plain}
\pagestyle{plain}
\pagenumbering{arabic}

\begin{abstract}
Consensus is fundamental for distributed systems since it underpins key functionalities of such systems ranging from distributed information fusion, decision-making, to decentralized control. In order to reach an agreement, existing consensus algorithms require each agent to exchange explicit state information with its neighbors. This leads to the disclosure of private state information, which is undesirable in cases where privacy is of concern. In this paper, we propose a novel approach for undirected networks which can enable secure and privacy-preserving average consensus in a decentralized architecture in the absence of an aggregator or third-party. By leveraging partial homomorphic cryptography to embed secrecy in pairwise interaction dynamics, our approach can guarantee convergence to the consensus value (subject to a quantization error) in a \emph{deterministic} manner without disclosing a node's state to its neighbors. We provide a new privacy definition for dynamical systems, and give a new framework to rigorously prove that a node's privacy can be protected as long as it has at least one legitimate neighbor which follows the consensus protocol faithfully without attempts to infer other nodes' states. In addition to enabling resilience to passive attackers aiming to steal state information, the approach also allows easy incorporation of defending mechanisms against active attackers who try to alter the content of exchanged messages. Furthermore, in contrast to existing noise-injection based privacy-preserving mechanisms which have to reconfigure the entire network when the topology or number of nodes varies, our approach is applicable to dynamic environments with time-varying coupling topologies. This secure and privacy-preserving approach is also applicable to weighted average consensus as well as maximum/minimum consensus under a new update rule. Numerical simulations and comparison with existing approaches confirm the theoretical results. Experimental results on a Raspberry-Pi board based micro-controller network are also presented to verify the effectiveness and efficiency of the approach.
\end{abstract}


\section{Introduction}
As a building block of distributed computing, average consensus as well as its various variants such as weighted average consensus and maximum/minimum consensus has been an active research topic in computer science and optimization for decades \cite{Morris74, Lynch96}. In recent years, with the advances of wireless communications and embedded systems, particularly the advent of wireless sensor networks and the Internet-of-Things, average consensus is finding increased applications in fields as diverse as automatic control, signal processing, social sciences, robotics, and optimization \cite{Olfati-Saber2007}.

Conventional consensus approaches employ the explicit exchange of state values among neighboring nodes to reach agreement on the computation. Such an explicit exchange of state information has two potential problems. First, it results in breaches of the privacy of participating nodes who may want to keep their state information confidential. For example, a group of individuals using average consensus to compute a common opinion may want to keep their personal opinions secret \cite{citeulike84}. Another example is power systems where multiple generators want to reach agreement on cost while keeping their individual generation information private \cite{Zhang11}. Secondly, storing or exchanging information in the unencrypted plaintext form is vulnerable to attackers who can steal information by hacking into the communication links or even the nodes. With the increased number of reported attack events and the growing awareness of security, keeping data encrypted in storage and communications has become the norm in many applications, particularly in many real-time sensing and control systems such as power systems and wireless sensor networks.

To address the pressing need for privacy and security in consensus, recently several solutions have been proposed. Most existing approaches use the idea of obfuscation to mask the true state values by adding noise on the state. Motivated by database privacy in computer science, \cite{Erfan15, Huang15, nozari2015, katewa2017privacy, huang2012differentially} exploited differential privacy to inject uncorrelated noise in average consensus. However, differential-privacy based approaches do not provide the exact average value due to their fundamental trade-off between enabled privacy and computational accuracy. To guarantee computational accuracy, \cite{Mo17, Manitara13, he2016private} proposed to inject correlated noise to exchanged messages to guarantee convergence to the exact value, which, however, also leads to vulnerabilities to external eavesdroppers, as shown in our numerical simulations in Fig. 8 and Fig. 9. Furthermore, these approaches normally rely on the assumption of a \textit{time-invariant} interaction graph, which is difficult to satisfy in many practical applications where the interaction patterns may change due to node mobility or fading communication channels. Observability based approaches have also been discussed to protect the privacy of multi-agent networks. The basic idea is to design the interaction topology so as to minimize the observability from a compromised agent, which amounts to minimizing its ability to infer the initial states of other network agents \cite{kia2015dynamic,observability1,observability2}. However, these approaches cannot protect the privacy of the direct neighbors of the compromised agent.

Neither can the aforementioned approaches protect nodes from active attackers who try to steal or even alter exchanged information by hacking into the nodes or the communication channels. To improve resilience to such attacks, a common approach is to employ cryptography. However, it is worth noting that although cryptography based approaches can easily provide privacy and security when an aggregator or third-party is available \cite{Gupta16}, like in cloud-based control or computation \cite{homomorphic_privacy1, homomorphic_privacy2, Lagendijk13}, their extension to completely \textit{decentralized} average consensus \textit{without any aggregators or third-parties} is extremely hard due to the difficulties in decentralized key management.

In this paper, we propose a homomorphic cryptography based approach for undirected networks that can guarantee the security and privacy of a node as long as it has at least one legitimate neighbor. Here, a legitimate neighbor is defined as a neighboring node who follows the consensus protocol faithfully without attempts to infer other nodes' states. Compared with differential-privacy based approaches in \cite{Erfan15, Huang15, nozari2015, katewa2017privacy, huang2012differentially} which trade accuracy for privacy and correlated-noise based approaches in \cite{Mo17, Manitara13, he2016private} which are vulnerable to external eavesdroppers, our approach can guarantee both the accuracy of consensus and resilience against external eavesdroppers. Despite the increased computational complexity caused by encryption, our approach is still manageable for resource constrained low-cost micro-controllers, as analyzed in detail in Sec. III and experimentally validated in Sec. VIII-C. Unlike the existing cryptography based average consensus approach in \cite{Lazzeretti14} and \cite{Freris16}, our approach allows every participating node to access the exact final value, and hence can be used in numerous applications where each agent has to use the precise consensus value to achieve cooperative control or information fusion.

The main contributions of the paper are as follows: 1) To our knowledge, our paper is the first to unify encryption and control dynamics in a completely decentralized manner without the assistance of a third party; 2) Our paper provides a new privacy definition for dynamical systems, and gives a new framework to rigorously prove that a node's privacy can be protected as long as it has at least one legitimate neighbor; 3) Our approach can be easily extended to provide security against active attackers aiming to alter the content of exchanged messages; and 4) We experimentally verified the efficiency of our approach using a Raspberry-Pi board based micro-controller network.

The outline of this paper is as follows. Sec. \ref{sec:background} reviews the average consensus problem and the homomorphic cryptography, particularly the Paillier cryptosystem. Our confidential interaction protocol is introduced in Sec. \ref{sec:proto}. In Sec. \ref{sec:analysis} we theoretically analyze the condition and rate of convergence, followed by a systematic discussion of privacy guarantees as well as security enforcement mechanisms in Sec. \ref{sec:security}. We further demonstrate in Sec. \ref{sec:application} that our confidential interaction protocol can be extended to weighted average consensus and maximum/minimum consensus under a new update rule. Some implementation issues are discussed in Sec. \ref{sec:impl}. Both numerical examples and hardware experiments (on Raspberry Pi boards) are presented in Sec. \ref{sec:example}. The conclusion is drawn in Sec. \ref{sec:conclusion}.

\section{Background}\label{sec:background}
In this section we briefly review the average consensus problem and the homomorphic encryption.

\subsection{Average Consensus}
We follow the same convention as in \cite{Olfati-Saber2007} where a network of $N$ agents is represented by an undirected graph ${G=(V,\,E,\,\mathbf{A})}$ with node set ${V}=\{v_1, \, v_2, \, \cdots v_N\}$, edge set ${E}\subset {V}\times {V}$, and a weighted adjacency matrix $\mathbf{A}=[a_{ij}]$ which satisfies $a_{ij}>0$ if $(v_i,v_j)\in E$ and $0$ otherwise. The set of neighbors of a node $v_i$ is denoted as
\begin{equation}
 {N}_i = \left\{v_j \in {V}| (v_i,v_j)\in {E}\right\}
\end{equation}
Throughout this paper we assume that the graph is undirected and connected. Therefore, $\mathbf{A}$ is symmetric
\begin{equation}\label{eq:aij}
 a^{(k)}_{ij} = a^{(k)}_{ji}>0\quad \forall (v_i,v_j)\in {E}
\end{equation}
Note that the superscript $k$ denotes that the weights are time-varying. Sometimes we drop $k$ for the sake of notational simplicity, but it is worth noting that all discussions in the paper are always applicable under time-varying weights. To achieve average consensus, namely converging of all states $x_i[k]$ $(i=1,2,\cdots,N)$ to the average of initial values, i.e., $\frac{\sum_{j=1}^N x_j[0]}{N}$, one commonly-used update rule is
\begin{equation}\label{eq:dt}
 x_i[k+1] = x_i[k] + \varepsilon\sum_{v_j\in N_i} a^{(k)}_{ij}(x_j[k]- x_i[k])
\end{equation}
where $\varepsilon$ is a constant step size.

\subsection{Homomorphic Encryption}
Our method to protect privacy and security is to encrypt the state. To this end, we briefly introduce a cryptosystem, more specifically the public-key cryptosystem which is applicable in open and dynamic networks without the assist of any trusted third party for key management. A public-key cryptosystem uses two keys -- a private key (also called secret key) and a public key distributed publicly. Any person can encrypt messages using a public key, but such messages can only be decrypted by the agents who have access to the private key. Most popular cryptosystems such as RSA \cite{Rivest1978}, ElGamal \cite{ElGamal1985}, and Paillier \cite{Paillier1999} are public-key cryptosystems. In this paper we focus on the Pailler cryptosystem. We adopt a simplified version of the Pailler cryptosystem from \cite{katz2008introduction} which has the following basic functions:
\begin{itemize}
 \item Key generation:
 \begin{enumerate}
 \item Choose two large prime numbers $p\in \mathbb{Z}$ and $q\in \mathbb{Z}$ of equal bit-length and compute $n=pq$.
 \item Let $\lambda= \phi(n)=(p-1)(q-1)$ where $\phi(\cdot)$ is the Euler's totient function.
 \item Let $\mu=\phi(n)^{-1}\;\text{mod}\;n$ which is the modular multiplicative inverse of $\phi(n)$.
 \item The public key $k_p$ is then $n$.
 \item The private key $k_s$ is then $(\lambda,\mu)$.
 \end{enumerate}
 \item Encryption ($c = \mathcal{E}(m)\in \mathbb{Z}_{n^2}^*$):\\
 Recall the definitions of $\mathbb{Z}_n= \{z|z \in \mathbb{Z}, 0 \leq z < n\}$ and $\mathbb{Z}^*_n = \{z|z \in \mathbb{Z}, 0 \leq z < n, \text{gcd}(z, n) = 1\}$ where $\text{gcd}(a,b)$ is the greatest common divisor of $a$ and $b$.
 \begin{enumerate}
 \item Choose a random $r\in \mathbb{Z}^*_n$.
 \item The ciphertext is given by $c=(n+1)^m  r^n\;\text{mod}\;n^2$, where $m\in \mathbb{Z}_n, c\in \mathbb{Z}^*_{n^2}$.
 \end{enumerate}
 \item Decryption ($m = \mathcal{D}(c)\in \mathbb{Z}_n$):
 \begin{enumerate}
 \item Define the integer division function $L(u) = \frac{u-1}{n}$.
 \item The plaintext is $m=L(c^\lambda\;\text{mod}\;n^2) \mu\;\text{mod}\;n$.
 \end{enumerate}
\end{itemize}

A cryptosystem is homomorphic if it allows certain computations to be carried out on the encrypted ciphertext. The Paillier cryptosystem is additive homomorphic because the ciphertext of $m_1+m_2$, i.e., $\mathcal{E}(m_1+m_2)$, can be obtained from $\mathcal{E}(m_1)$ and $\mathcal{E}(m_2)$ directly:
\begin{equation} \label{eq:add_explicity}
 \begin{aligned}
 \mathcal{E}(m_1,r_1) \mathcal{E}(m_2,r_2) =& ((n+1)^{m_1} {r_1}^n) ((n+1)^{m_2} {r_2}^n)\;\text{mod}\;n^2\\
 =&((n+1)^{m_1+m_2}(r_1r_2)^n)\;\text{mod}\;n^2\\
 =&\mathcal{E}(m_1+m_2, r_1r_2)
 \end{aligned}
\end{equation}
The dependency on random numbers $r_1$ and $r_2$ is explicitly shown in (\ref{eq:add_explicity}), yet they play no role in the decryption. Therefore, the following shorthand notation will be used instead:
\begin{equation}\label{eq:add}
 \mathcal{E}(m_1) \mathcal{E}(m_2)=\mathcal{E}(m_1+m_2)
\end{equation}
Moreover, if we multiply the same ciphertext $k\in\mathbb{Z}^+$ times, we can obtain
\begin{equation}\label{eq:mult}
 \mathcal{E}(m)^k = \prod_{i=1}^k \mathcal{E}(m)=\mathcal{E}(\sum_{i=1}^km)=\mathcal{E}(km)
\end{equation}
Notice however, the Paillier cryptosystem is not multiplicative homomorphic because $k$ in (\ref{eq:mult}) is in the plaintext form. Furthermore, the existence of the random number $r$ in Paillier cryptosystem gives it resistance to dictionary attacks \cite{Goldreich_2} which infer a key to an encrypted message by systematically trying all possibilities, like exhausting all words in a dictionary. Moreover, since Paillier cryptography only works on numbers that can be represented by binary strings, we multiply a real-valued state by a large integer $Q$ before converting it to a binary string so as to ensure small quantization errors. The details will be discussed in Sec. VII-A.

\section{Confidential Interaction Protocol}\label{sec:proto}
In this section, we propose a completely decentralized, third-party free confidential interaction protocol that can guarantee average consensus while protecting the privacy of all participating nodes. Instead of adding noise to hide the states, our approach combines encryption with randomness in the system dynamics, i.e., the coupling weights $a^{(k)}_{ij} \in \mathbb{R}$, to prevent two communicating parties in a pairwise interaction from exposing information to each other. In this way the states are free from being contaminated by covering noise, guaranteeing a deterministic convergence to the average (subject to a quantization error). The computational complexity of this protocol is also discussed in this section.

We first present details of our confidential interaction protocol based on (\ref{eq:dt}). In particular we show how a node can obtain the weighted difference (\ref{eq:diff}) between itself and any of its neighbor(s) without disclosing each other's state information:
\begin{equation} \label{eq:diff}
    \begin{aligned}
 \Delta x_{ij}[k] =& a_{ij}^{(k)} (x_j[k]-x_i[k])\\
 \Delta x_{ji}[k] =& a_{ji}^{(k)} (x_i[k]-x_j[k])\\
 \text{subject to} & \quad a_{ij}^{(k)} = a_{ji}^{(k)}>0
 \end{aligned}
\end{equation}
Plugging the state difference (\ref{eq:diff}) into (\ref{eq:dt}) gives a new formulation of average consensus
\begin{equation} \label{eq:ex_dt}
 x_i[k+1] = x_i[k] + \varepsilon\sum_{v_j\in N_i} \Delta x_{ij}[k]
\end{equation}
Notice that in a decentralized system it is impossible to protect the privacy of both nodes in a pairwise interaction if the protocol (\ref{eq:diff}) is used without a third party distributing secret $a^{(k)}_{ij}$. This is due to the fact that even if we encrypt all the intermediate steps, if one node, for instance $v_i$, has access to $a^{(k)}_{ij}$, it can still infer the value of $x_j[k]$ through $x_j[k] = {\Delta x_{ij}[k]}/{a^{(k)}_{ij}}+x_i[k]$. From now on, for the sake of simplicity in bookkeeping, we omit the superscript $k$ in $a_{ij}^{(k)}$. But it is worth noting that all the results hold for time-varying weights.

\begin{figure}
 \centering
 \includegraphics[width=.48\textwidth]{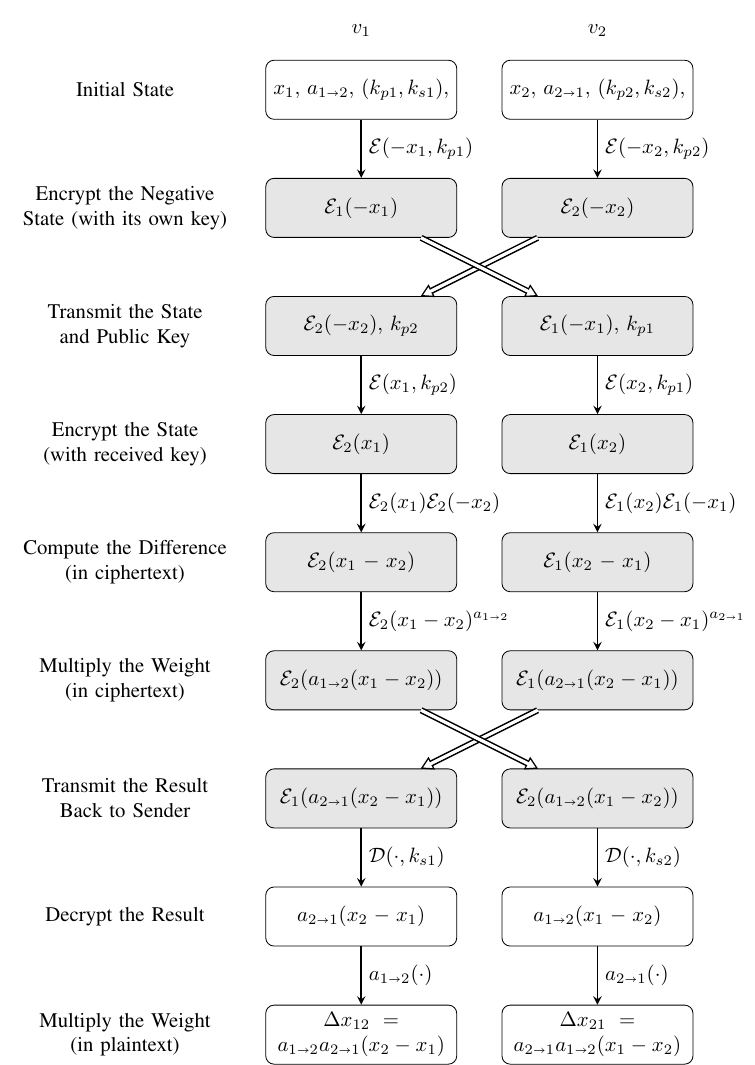}
 \caption{A step-by-step illustration of the confidential interaction protocol. Single arrows indicate the flow of computations; double arrows indicate data exchange via a communication channel. Shaded nodes indicate the computation done in ciphertext. Note that $a_{1 \shortrightarrow 2}$ and $a_{2 \shortrightarrow 1}$ are different from step to step.}
 \label{fig:1}
\end{figure}

We solve this problem by constructing each weight $a_{ij}$ as the product of two random numbers, namely $a_{ij}=a_{ji}=a_{i \shortrightarrow j}  a_{j \shortrightarrow i}$ with $a_{i \shortrightarrow j}$ randomly generated by and only known to node $v_i$ and $a_{j \shortrightarrow i}$ randomly generated by and only known to node $v_j$. We will show later that this weight construction approach renders two interacting nodes unable to infer each other's state while guaranteeing convergence to the average. Next, without loss of generality, we consider a pair of connected nodes ($v_1,\,v_2$) to illustrate the idea (cf. Fig. \ref{fig:1}). For simplicity, we assume that the states $x_1$ and $x_2$ are scalar. Each node maintains its own public and private key pairs $(k_{pi}, k_{si}),\; i\in\{1,\,2\}$.

Due to symmetry, we only show how node $v_1$ obtains the weighted state difference, i.e., the flow $v_1 \rightarrow v_2 \rightarrow v_1$. Before starting the information exchange, node $v_1$ (resp. $v_2$) generates its new random number $a_{1 \shortrightarrow 2}$ (resp. $a_{2 \shortrightarrow 1}$). First, node $v_1$ sends its encrypted negative state $\mathcal{E}_1(-x_1)$ as well as the public key $k_{p1}$ to node $v_2$. Note that here the subscript in $\mathcal{E}_1$ denotes encryption using the public key of node $v_1$. Node $v_2$ then computes the encrypted $a_{2 \shortrightarrow 1}$-weighted difference $\mathcal{E}_1\left(a_{2 \shortrightarrow 1}(x_2-x_1)\right)$ following the three steps below:
\begin{enumerate}
 \item Encrypt $x_2$ with $v_1$'s public key $k_{p1}$: $x_2 \rightarrow \mathcal{E}_1(x_2)$.
 \item Compute the difference directly in ciphertext:
 \begin{equation}
 \mathcal{E}_1(x_2-x_1)=\mathcal{E}_1(x_2+(-x_1))=\mathcal{E}_1(x_2)  \mathcal{E}_1(-x_1)
 \end{equation}
 \item Compute the $a_{2 \shortrightarrow 1}$-weighted difference in ciphertext:
 \begin{equation}
 \label{eq:exp_a1}
 \mathcal{E}_1\left(a_{2 \shortrightarrow 1}  (x_2-x_1)\right)=\left(\mathcal{E}_1(x_2-x_1)\right)^{a_{2 \shortrightarrow 1}}
 \end{equation}
\end{enumerate}
Then $v_2$ returns $\mathcal{E}_1\left(a_{2 \shortrightarrow 1}  (x_2-x_1) \right)$ to $v_1$. After receiving $ \mathcal{E}_1 \left(a_{2 \shortrightarrow 1}   (x_2-x_1)\right)$, $v_1$ decrypts it using the private key $k_{s1}$ and multiplies the result with $a_{1 \shortrightarrow 2}$ to get the weighted difference $\Delta x_{12}$:
\begin{equation}\label{eq:mul_a1}
 \begin{aligned}
 \mathcal{E}_1\left(a_{2 \shortrightarrow 1} (x_2-x_1)\right) \xrightarrow {\mathcal{D}_1} a_{2 \shortrightarrow 1}(x_2-x_1)\\
 \Delta x_{12}=a_{1 \shortrightarrow 2}a_{2 \shortrightarrow 1}(x_2-x_1)
 \end{aligned}
\end{equation}

In a similar manner, the exchange $v_2\rightarrow v_1 \rightarrow v_2$ produces $\mathcal{E}_2 \left(a_{1 \shortrightarrow 2}(x_1-x_2)\right)$ for $v_2$ who then decrypts the message and multiplies the result by its own multiplier $a_{2 \shortrightarrow 1}$ to get $\Delta x_{21}$
\begin{equation} \label{eq:mul_a2}
 \begin{aligned}
 \mathcal{E}_2 \left(a_{1 \shortrightarrow 2} (x_1-x_2)\right) \xrightarrow{\mathcal{D}_2}a_{1 \shortrightarrow 2}(x_1-x_2)\\
 \Delta x_{21}=a_{2 \shortrightarrow 1} a_{1 \shortrightarrow 2}(x_1-x_2)
 \end{aligned}
\end{equation}
After each node collects the weighted differences from all neighbors, it updates its state with (\ref{eq:ex_dt}) accordingly.

Several remarks are in order:
\begin{itemize}
 \item The construction of each $a_{ij}$ as the product of two random numbers $a_{i \shortrightarrow j}$ and $a_{j \shortrightarrow i}$ is key to guarantee that the weights are symmetric, i.e., $a_{ij}=a_{ji}$, which is necessary for average consensus. Note that the coupling weights $a_{i \shortrightarrow j}$ are randomly chosen by individual nodes and are completely independent of communication channels.
 \item In the initial iteration $k=0$, node $v_i \in V$ chooses $a_{i \shortrightarrow j}$ randomly from the set of real numbers $\mathbb{R}$, which is critical for privacy preservation, as will be clear in the proof of Theorem 3. After the initial iteration, i.e., for $k \geq 1$, node $v_i \in V$ chooses $a_{i \shortrightarrow j}$ randomly from a finite admissible range $[\barbelow{a}, \, \bar{a}]$. Note that $a_{i \shortrightarrow j}$ can have any random distribution in $[\barbelow{a}, \, \bar{a}]$ such as uniform, reciprocal, and truncated Gaussian.
 \item To guarantee the convergence of average consensus, the admissible range $[\barbelow{a}, \, \bar{a}]$ must satisfy $0 < \barbelow{a} < \bar{a} < \frac{1}{\sqrt{\varepsilon \Delta}}$ where $\varepsilon$ is the constant step size in (\ref{eq:dt}) and $\Delta \triangleq \max_{i}|N_i|$ with $|\bullet|$ denoting the set cardinality. Therefore, a smaller $\varepsilon$ means a larger $\bar{a}$, which results in an extended admissible range $[\barbelow{a}, \, \bar{a}]$ in the selection of the random number $a_{i \shortrightarrow j}$.
 \item $v_2$ does not have the private key of $v_1$ and cannot see $x_1$ which is encrypted in $\mathcal{E}_1(-x_1)$.
 \item Given $a_{2 \shortrightarrow 1}(x_2-x_1)$, $v_1$ cannot solve for $x_2$ because $a_{2 \shortrightarrow 1}$ is only known to $v_2$.
 \item At each iteration, real-valued states are converted to fixed point representation for encryption; the weighted differences are converted back to real values for update.
 \item We encrypt $\mathcal{E}_1(-x_1)$ because it is much more difficult to compute subtraction in ciphertext. The issue regarding encrypting negative values using Paillier is discussed in Sec. \ref{sec:impl}.
 \item Different from \cite{Lazzeretti14} and \cite{Freris16} which deprive individual nodes from accessing the average value (note that in \cite{Freris16} individual participating nodes do not have the decryption key to decrypt the final consensus value which is in the encrypted form, otherwise they will be able to decrypt intermediate computations to access other nodes' states), our approach enables every node to calculate the average value in a completely decentralized manner, and hence can be used in numerous applications where each agent relies on the precise consensus value to achieve cooperative control or information fusion.
\end{itemize}

Next we analyze the complexity of our confidential interaction protocol. Note that, although the key generation process is typically the most time consuming step, it incurs a one-time fixed cost. Unless otherwise stated, we assume that the keys are reused for subsequent steps.

\begin{Lemma 1}
    Under our confidential interaction protocol illustrated in Fig. \ref{fig:1}, for a particular node $v_i$ in a connected network, the total computation overhead of each iteration is $\mathcal{O}(|N_i| l)$ where $|N_i|$ is the number of neighboring nodes of node $v_i$ and $l$ is the bit length of the public key $n$.
\end{Lemma 1}

\textit{Proof}: The encryption step requires two modulo exponentiation steps which can be computed in $\mathcal{O}(l+\log(m))$ time, where $m$ is the message's integer value and $l$ is the bit length of the public key $n$ \cite{goldreich2007foundations}. The multiplication and modulo steps are both $\mathcal{O}(1)$. In practice the length of the public key $n$ is much longer than the size of $m$, so the complexity of the encryption is assumed to be $\mathcal{O}(l)$.

The decryption step requires one modulo exponentiation, one subtraction, one division followed by a multiplication and a modulo operation. The overall time complexity is $\mathcal{O}(l)$.

During each message exchange, node $v_i$ must encrypt its own message once and decrypt $|N_i|$ return messages, so the overall complexity is $\mathcal{O}((|N_i|+1)l)$. For each $|N_i|$ in-bound messages, node $v_i$ has to encrypt its own message, compute the difference, and then modulo-exponentiate the difference. Therefore the total computation complexity is $\mathcal{O}(|N_i|(l + \log(\bar{a})))$, where $\bar{a}$ is the upper bound of the admissible range for the random selection of $a_{i \shortrightarrow j}^{(k)}$.

As a result, for a particular node $v_i$, the total computation overhead of each iteration is $\mathcal{O}(|N_i| l)$.
\hfill{$\blacksquare$}

From Lemma 1, it can be seen that the computational complexity of our approach does not increase with network size but rather with the number of neighbors (so the computational complexity of our approach is moderate even for large networks with moderate connections). In Sec. VIII-C, we will use a network of Raspberry Pi boards to experimentally verify that the computational complexity is easily manageable on resource-constrained micro-controller based real-time control systems.

\section{Theoretical Analysis of Convergence}\label{sec:analysis}

In this section, we first show that average consensus will be guaranteed under the confidential interaction protocol. Then we provide convergence speed analysis of our approach.

Let $\mathbf{x}\in \mathbb{R}^n$ denote the augmented state vector of all nodes. The network dynamics in (\ref{eq:dt}) can be rewritten as:
\begin{equation} \label{eq:perron}
 \mathbf{x}[k+1]=\mathbf{P}^{(k)}\mathbf{x}[k]
\end{equation}
where $\mathbf{P}^{(k)}=\mathbf{I}-\varepsilon \mathbf{L}^{(k)}$ is the Perron matrix and $\mathbf{L}^{(k)}=[l^{(k)}_{ij}]$ is the time-varying Laplacian matrix defined by
\begin{equation}\label{eq:weithts in DT}
\begin{aligned}
    l^{(k)}_{ij}=&
    \begin{cases}
        \sum_{j\in N_i} a^{(k)}_{ij} & i= j\\ -a^{(k)}_{ij} & i \neq j
    \end{cases}
\end{aligned}
\end{equation}

\begin{Theorem 1}
For a connected network of $N$ nodes, if the coupling weights $a^{(k)}_{ij}$ in (\ref{eq:weithts in DT}) are established according to the confidential interaction protocol illustrated in Fig. \ref{fig:1} and the admissible range $[\barbelow{a}, \, \bar{a}]$ for the random selection of $a_{i \shortrightarrow j}^{(k)}$ satisfies $0 < \barbelow{a} < \bar{a} < \frac{1}{\sqrt{\varepsilon \Delta}}$, then the system will achieve average consensus with states converging to
\begin{equation}\label{consensus_value}
\lim_{k\to\infty}\mathbf{x}[k] = \alpha \mathbf{1}\quad\text{with }\alpha=\text{Avg}[0]=\frac{1}{N}\mathbf{1}^T\mathbf{x}[0]
\end{equation}
\end{Theorem 1}

\textit{Proof}: The proof can be obtained by following the similar line of reasoning of Theorem 2 in \cite{Olfati-Saber2007}.
\hfill{$\blacksquare$}

\begin{Remark 1}
 Since the framework allows time-varying weighted matrix $\mathbf{A}^{(k)}$, it can easily be extended to the case with switching interaction graphs according to \cite{Moreau05}.
\end{Remark 1}

To analyze the convergence speed of our approach, following the convergence definition in \cite{nedic2017achieving} and \cite{nedich2016geometrically}, we define the convergence rate as at least $\gamma \in (0, \, 1)$ if there exists a positive constant $C$ such that $\big\|\mathbf{x}[k] - \alpha \mathbf{1} \big\| \leq C \gamma ^k$ holds for all $k$, where $\alpha$ is the average value of the initial states of all nodes in (\ref{consensus_value}). Note that a smaller $\gamma$ means a faster convergence speed.

\begin{Theorem 2}
    For a connected network of $N$ nodes, if the coupling weights $a^{(k)}_{ij}$ in (\ref{eq:weithts in DT}) are established according to the confidential interaction protocol illustrated in Fig. \ref{fig:1} and the admissible range $[\barbelow{a}, \, \bar{a}]$ for the random selection of $a_{i \shortrightarrow j}^{(k)}$ satisfies $0 < \barbelow{a} < \bar{a} < \frac{1}{\sqrt{\varepsilon \Delta}}$, then the convergence rate of our privacy-preserving average consensus algorithm is at least $\gamma= (1-\eta^{N-1}) ^{\frac{1}{N-1}} \in (0, \, 1)$ with $\eta=\min \{1- \varepsilon \Delta \bar{a}^2, \, \varepsilon \barbelow{a}^2 \}$, meaning that there exists a positive constant $C$ such that $\big\|\mathbf{x}[k] - \alpha \mathbf{1} \big\| \leq C \gamma ^k$ holds for all $k$.
\end{Theorem 2}

{\it Proof}: For iteration $k=0$, since $\mathbf{1}^T \mathbf{P}^{(0)} = \mathbf{1}^T$ holds, we have $\mathbf{1}^T \mathbf{x}[1]= \mathbf{1}^T \mathbf{P}^{(0)} \mathbf{x}[0] = \mathbf{1}^T \mathbf{x}[0]$. So we have $\alpha = \frac{1}{N} \sum_{j=1}^{N} {x}_j[0] = \frac{1}{N} \sum_{j=1}^{N} {x}_j[1]$.

For iteration $k \geq 1$, each random number $a_{i \shortrightarrow j}^{(k)}$ is randomly chosen from $[\barbelow{a}, \, \bar{a}]$ by node $i$ with $0< \barbelow{a}< \bar{a} < \frac{1}{\sqrt{\varepsilon \Delta}}$. Further taking into account $\mathbf{P}^{(k)}=\mathbf{I}-\varepsilon \mathbf{L}^{(k)}$, we have the following properties for $\mathbf{P}^{(k)}$'s $(i,i)$th element $p_{ii}^{(k)}$ and $(i,j)$th element $p_{ij}^{(k)}$ for any $1\leq i,j \leq N$:
\begin{equation}\label{p_ij_1}
p_{ii}^{(k)} =1- \varepsilon \sum\limits_{j \in N_i} a_{i \shortrightarrow j}^{(k)} a_{j \shortrightarrow i}^{(k)} \geq 1- \varepsilon \Delta \bar{a}^2 \qquad \forall \, v_i \in V
\end{equation}
and
\begin{equation}\label{p_ij_2}
\left\{\begin{aligned}
p_{ij}^{(k)}& = 0 \qquad \qquad \qquad \qquad \forall \, v_i \in V \, \text{and} \, v_j \notin N_i \cup \{v_i\}\\
p_{ij}^{(k)}& = \varepsilon a_{i \shortrightarrow j}^{(k)} a_{j \shortrightarrow i}^{(k)} \geq \varepsilon \barbelow{a}^2 \quad \ \, \forall \, v_i \in V \, \text{and} \, v_j \in N_i
\end{aligned}\right.
\end{equation}

Defining $\eta$ as $\eta \triangleq \min \{1- \varepsilon \Delta \bar{a}^2, \, \varepsilon \barbelow{a}^2 \}$, we have $0<\eta<1$ due to $0< \barbelow{a}< \bar{a} < \frac{1}{\sqrt{\varepsilon \Delta}}$. According to (\ref{p_ij_1}) and (\ref{p_ij_2}), we have $p_{ij}^{(k)} \geq \eta$ for $v_i \in V$ and $v_j \in N_i \cup \{v_i\}$.

Denote $\mathbf{P}^{(k)} \mathbf{P}^{(k-1)} \cdots \mathbf{P}^{(s)}$ as $\mathbf{\Phi}(k,s)$. Note that $\mathbf{\Phi}(k,k)= \mathbf{P}^{(k)}$ holds for all $k$. Since $\mathbf{P}^{(k)}$ ($k \geq 1$) satisfies the Assumptions 2, 3, 4, and 5 in \cite{nedic2010constrained}, following its Proposition 1, we have
\begin{equation}\label{phi_i_j}
\begin{aligned}
\Big|{\phi}^{(k,s)}_{ij}-\frac{1}{N}\Big| \leq 2 \frac{1+\eta^{-N+1}}{1-\eta^{N-1}}(1-\eta^{N-1})^{\frac{k-s}{N-1}}
\end{aligned}
\end{equation}
for $i, j=1,2,\ldots,N$ and $k \geq s \geq 1$ where ${\phi}^{(k,s)}_{ij}$ represents the $(i,j)$th element of $\mathbf{\Phi}(k,s)$.

From (\ref{eq:perron}), we have
\begin{equation}\label{x_k_1}
\begin{aligned}
\mathbf{x}[k] = \mathbf{P}^{(k-1)} \mathbf{P}^{(k-2)} \cdots \mathbf{P}^{(1)} \mathbf{x}[1] = \mathbf{\Phi}(k-1,1) \mathbf{x}[1]
\end{aligned}
\end{equation}
for $k \geq 2$. Combining (\ref{p_ij_1}) and (\ref{p_ij_2}) with (\ref{x_k_1}), we obtain
\begin{equation}
\begin{aligned}
\big\| \mathbf{x}[k] - \alpha \mathbf{1} \big\| & = \big\| \mathbf{\Phi}(k-1,1) \mathbf{x}[1] - \alpha \mathbf{1} \big\| \\
& = \Big( \sum\limits_{i=1}^{N} \Big| \sum\limits_{j=1}^{N} ( \phi^{(k-1,1)}_{ij} - \frac{1}{N}) {x}_j[1] \Big|^2 \Big)^{1/2}\\
& \leq \Big( \sum\limits_{i=1}^{N} \big( \sum\limits_{j=1}^{N} \Big|\phi^{(k-1,1)}_{ij} - \frac{1}{N}\Big|^2 \big\| \mathbf{x}[1] \big\|^2 \big) \Big)^{1/2}\\
& \leq 2N\frac{1+\eta^{-N+1}}{1-\eta^{N-1}}(1-\eta^{N-1})^{\frac{k-2}{N-1}} \big\| \mathbf{x}[1] \big\|
\end{aligned}
\end{equation}
for $k \geq 2$, where we used Cauchy-Schwarz inequality in the derivation. By re-arranging the terms, we have
\begin{equation}
\begin{aligned}
\big\|\mathbf{x}[k] - \alpha \mathbf{1} \big\| \leq C_0 \gamma ^k
\end{aligned}
\end{equation}
for $k \geq 2$, where
\begin{equation}
\begin{aligned}
C_0 = 2 N \big\| \mathbf{x}[1] \big\| \frac{1+\eta^{-N+1}} {(1-\eta^{N-1})^{\frac{N+1}{N-1}}}
\end{aligned}
\end{equation}
and
\begin{equation}\label{gamma_def}
\begin{aligned}
\gamma= (1-\eta^{N-1}) ^{\frac{1}{N-1}}
\end{aligned}
\end{equation}
Given $\eta\in (0, \, 1)$, it can be easily verified that $\gamma \in (0,\,1)$ holds.

Defining $C$ as
\[C \triangleq \max\{ \big\|\mathbf{x}[0]- \alpha \mathbf{1}\big\|, \, \big\|\mathbf{x}[1]- \alpha \mathbf{1}\big\| \gamma^{-1}, \, C_0\}
\]
we have
\begin{equation}\label{convergence_def}
\begin{aligned}
\big\|\mathbf{x}[k] - \alpha \mathbf{1} \big\| \leq C \gamma ^k
\end{aligned}
\end{equation}
for all $k$. Therefore, the convergence rate of our privacy-preserving average consensus algorithm is at least $\gamma= (1-\eta^{N-1}) ^{\frac{1}{N-1}} \in (0, \, 1)$.
\hfill{$\blacksquare$}

From Theorem 2, it can be seen that the initial weights will affect the absolute convergence time as they affect the parameter $C$ in (\ref{convergence_def}). However, they will not affect the exponential convergence rate $\gamma$ according to (\ref{gamma_def}). According to (\ref{convergence_def}), it is easy to see that a smaller $\gamma$ leads to a faster convergence. Given the relationship $\gamma= (1-\eta^{N-1}) ^{\frac{1}{N-1}}$, to get a faster convergence speed, i.e., a smaller $\gamma$, it suffices to increase $\eta$. Further noting $\eta=\min \{1- \varepsilon \Delta \bar{a}^2, \, \varepsilon \barbelow{a}^2 \}$, increasing $\eta$ amounts to decreasing $\bar{a}$ and/or increasing $\barbelow{a}$, which in turn reduces the size of the range $[\barbelow{a}, \, \bar{a}]$. Therefore, there is a trade-off between the size of the admissible range $[\barbelow{a}, \, \bar{a}]$ and the convergence speed. However, this trade-off only matters if one wants to protect the intermediate states from being inferrable by honest-but-curious attackers (a reduced $[\barbelow{a}, \, \bar{a}]$ reduces the admissible range for randomizing $a_{i \shortrightarrow j}$ and $a_{j \shortrightarrow i}$, and hence enables an honest-but-curious attacker $v_i$ to get a better range estimation of its neighbor $v_j$'s intermediate states $x_j[k]$ ($k \geq 1$) based on received $\Delta x_{ij}$). If the intermediate states do not need to be protected as in the consensus problem considered in this paper (where only the initial state values matter), we can still enhance convergence speed by reducing the size of the admissible range $[\barbelow{a}, \, \bar{a}]$. In fact, even if we reduce the size of $[\barbelow{a}, \, \bar{a}]$ to zero for iterations after $k=0$, i.e., setting $\barbelow{a}=\bar{a}$ starting from $k=1$ (note that in the initial iteration $k=0$ coupling weights $a^{(0)}_{i \shortrightarrow j}$ should be randomly chosen from $\mathbb{R}$), our approach can still guarantee the privacy of initial state values, as is clear from the proof of Theorem 3 in Sec. V-A.

\section{Analysis of Privacy and Security}\label{sec:security}
Privacy and security are sometimes used interchangeably in the literature but here we make the distinction explicit. Among the control community privacy is equivalent to the concept of unobservability. Privacy is also closely related to the concept of semantic security from cryptography \cite{Goldreich_2}. Both concepts essentially concern with an honest-but-curious adversary which is interested in learning the states of the network but conforms to the rules of the system. Security, on the other hand, deals with a broader issue which includes learning the states as well as the possibilities of exploiting the system to cause damages.

\subsection{Privacy Guarantees}

Our protocol provides protection against an honest-but-curious adversary, which can be a node in the network or an observer eavesdropping communication links.

The Paillier encryption algorithm is known to provide semantic security, i.e., Indistinguishability under Chosen Plaintext Attack (IND-CPA) \cite{Paillier1999}. As a result, the recipient of the first transmission $\mathcal{E}(-x_i)$ cannot see the value of $x_i$ at all time.

Before analyzing the privacy-preserving performance of our approach, we first give a definition of privacy preservation used throughout this paper.

\begin{Definition 1}\label{def_1}
For a connected network of $N$ nodes, the privacy of the initial value $x_i[0]$ of node $v_i$ is preserved if an honest-but-curious adversary cannot estimate the value of $x_i[0]$ with any accuracy.
\end{Definition 1}

Definition \ref{def_1} requires that an honest-but-curious adversary cannot even find a range for a private value and thus is more stringent than the privacy preservation definition considered in \cite{liu2006random, han2010privacy, cao2014privacy} which defines privacy preservation as the inability of an adversary to {\it uniquely} determine the protected value. It is worth noting that by a finite range, we mean lower and upper bounds that an adversary may obtain based on accessible information. We do not consider representation bounds caused by the finite number of bytes that can be used to represent a number in a computer.

As per the naming convention in cryptography, it is customary to name the legitimate sender and receiver participants as $A$ (Alice) and $B$ (Bob), and the adversary as $E$ (Eve).

\begin{figure}
    \centering
    \includegraphics[width=0.48\textwidth]{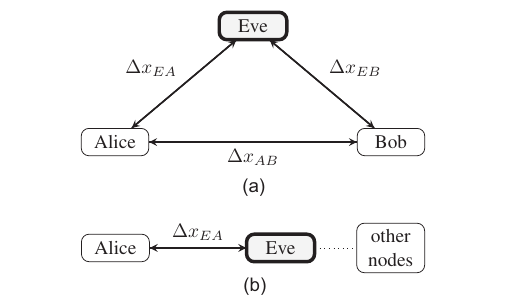}
    \caption{Two connection configurations considered in the proofs of Theorem 3 (subplot a) and Theorem 4 (subplot b).}
    \label{fig:privacy}
\end{figure}

\begin{Theorem 3}\label{thm:3}
 For a connected network of $N$ nodes with system dynamics in (\ref{eq:dt}), we assume that all nodes follow the confidential interaction protocol illustrated in Fig. \ref{fig:1}. An honest-but-curious node Eve who can receive messages from a neighboring node Alice cannot learn the initial state of Alice if Alice is also connected to another legitimate node Bob.
\end{Theorem 3}

\textit{Proof}: Without loss of generality, we consider the connection configuration illustrated in Fig. \ref{fig:privacy}(a) where Eve can interact with both Alice and Bob. If Eve cannot infer the state of Alice or Bob in this configuration, neither can it when either the Alice--Eve connection or the Bob--Eve connection is removed which reduces the amount of information accessible to Eve.

We propose a new privacy-proving approach based on the indistinguishability of a private value's arbitrary variation to Eve. We define the information accessible to Eve at iteration $k$ as $I_E[k]=\big\{ \Delta x_{EA}[k], \, \Delta x_{EB}[k], \, x_E[k], \, a^{(k)}_{E \shortrightarrow A}, \, a^{(k)}_{E \shortrightarrow B} \big\}$. So as time evolves, the cumulated information accessible to Eve can be summarized as $I_E = \bigcup_{k=0}^{\infty} I_E[k]$.

To show that the privacy of the initial value $x_A[0]$ can be preserved against Eve, i.e., Eve cannot estimate the value of $x_A[0]$ with any accuracy, it suffices to show that under any initial value $\bar{x}_A[0] \neq x_A[0]$ the information accessible to Eve, i.e., $\bar{I}_E$ could be exactly the same as $I_E$, the cumulated information accessible to Eve under $x_A[0]$. This is because the only information available for Eve to infer the initial value $x_A[0]$ is $I_E$, and if $I_E$ could be the outcome under any initial values of $x_A[0]$, then Eve has no way to even find a range for the initial value $x_A[0]$. Therefore, we only need to prove that for any $\bar{x}_A[0] \neq x_A[0]$, $\bar{I}_E = I_E$ could hold.

Next we show that there exist initial values of $x_B[0]$ and coupling weights satisfying the requirements of the confidential interaction protocol in Sec. III that make $\bar{I}_E = I_E$ hold under $\bar{x}_A[0] \neq x_A[0]$. (Note that the alternative initial values of $x_B[0]$ should guarantee that the agents still converge to the original average value after $x_A[0]$ is altered to $\bar{x}_A[0]$.) More specifically, under the following initial condition
\begin{equation}\label{initial_values}
\begin{aligned}
\bar{x}_B[0] = x_A[0] + x_B[0] -\bar{x}_A[0]
\end{aligned}
\end{equation}
and coupling weights
\begin{equation}\label{coupling_weight_bar}
\begin{aligned}
& \bar{a}^{(0)}_{E \shortrightarrow A}=a^{(0)}_{E \shortrightarrow A} \\
& \bar{a}^{(0)}_{E \shortrightarrow B}=a^{(0)}_{E \shortrightarrow B} \\
& \bar{a}^{(0)}_{A \shortrightarrow E}=a^{(0)}_{A \shortrightarrow E}(x_E[0]-x_A[0])/(x_E[0]-\bar{x}_A[0]) \\
& \bar{a}^{(0)}_{A \shortrightarrow B}= (x_A[0] - \bar{x}_A[0]) / ( \varepsilon a^{(0)}_{B \shortrightarrow A} (x_B[0] - x_A[0]) ) + a^{(0)}_{A \shortrightarrow B}\\
& \bar{a}^{(0)}_{B \shortrightarrow E}=a^{(0)}_{B \shortrightarrow E}(x_E[0]-x_B[0])/(x_E[0]- \bar{x}_B[0]) \\
& \bar{a}^{(0)}_{B \shortrightarrow A}=a^{(0)}_{B \shortrightarrow A}(x_B[0]-x_A[0]) / (\bar{x}_B[0]-\bar{x}_A[0]) \\
& \bar{a}^{(k)}_{i \shortrightarrow j}=a^{(k)}_{i \shortrightarrow j} \quad \forall \, i, j \in \{A, B, E \}, \ k=1, 2, \ldots\\
\end{aligned}
\end{equation}
it can be easily verified that $\bar{I}_E = I_E$ holds for any $\bar{x}_A[0] \neq x_A[0]$. Note that (\ref{initial_values}) is used to guarantee that the consensus value does not change under the new initial values $\bar{x}_A[0]$ and $\bar{x}_B[0]$. Therefore, the honest-but-curious node Eve cannot learn the initial state of a neighboring node Alice based on accessible information if Alice is also connected to another legitimate node Bob. \hfill{$\blacksquare$}

It is worth noting that conventional observability based approaches in e.g., \cite{observability1} and \cite{observability2} cannot be used to analyze node Eve's ability to infer the states of Alice because the coupling weights $a^{(k)}_{A \shortrightarrow E}$ are not available to Eve.

\begin{Remark 2}
In the confidential interaction protocol, we allow the coupling weights for $k=0$ to be unrestricted and random chosen from $\mathbb{R}$. This unrestricted admissible range for coupling weights of $k=0$ is key to guarantee that no matter what value $\bar{x}_A[0]$ is, the resulting $\bar{a}^{(0)}_{E \shortrightarrow A}$, $\bar{a}^{(0)}_{E \shortrightarrow B}$, $\bar{a}^{(0)}_{A \shortrightarrow E}$, $\bar{a}^{(0)}_{A \shortrightarrow B}$, $\bar{a}^{(0)}_{B \shortrightarrow E}$, and $\bar{a}^{(0)}_{B \shortrightarrow A}$ in the left-hand-side of (\ref{coupling_weight_bar}) are always admissible and hence guarantee the achievement of the defined privacy.
\end{Remark 2}

\begin{Remark 3}
From the derivation in Theorem 3, one can get that if only initial states need to be protected, as defined in Definition 1, then encrypting messages for the initial iteration $k=0$ is sufficient. Avoiding encrypting messages after the initial iteration can greatly reduce the computational complexity to a one-time cost and is highly desirable when computational complexity is of concern. Of course, if all intermediate messages need to be protected from being disclosed to an external eavesdropper, we still need to encrypt exchanged messages after the initial iteration.
\end{Remark 3}

\begin{Remark 3}
Given that the random distribution of $a_{i \shortrightarrow j}$ can be non-stationary, the approach is also resilient to statistical inferences based on law-of-large-numbers \cite{konheim2007computer}.
\end{Remark 3}

\begin{Remark 4}
From the derivation, we can see that when Alice is connected to at least one legitimate node Bob, the collusion of multiple honest-but-curious nodes does not help the inference of Alice's initial state.
\end{Remark 4}

\begin{Remark 5}
Following the same line of reasoning, it can be obtained that an honest-but-curious node Eve 1 cannot infer the initial state of a neighboring node Alice if Alice is also connected to another honest-but-curious node Eve 2 that does not collude with Eve 1.
\end{Remark 5}

\begin{Remark 6}
    Our approach can enable privacy even when Eve can interact with all neighbors of Alice, which is not allowed in existing approaches in \cite{Mo17} and \cite{Manitara13}.
\end{Remark 6}

Based on the analysis framework, we can also obtain a situation where it is possible for Eve to infer other nodes' states which should be avoided.

\begin{Theorem 4}\label{thm:4}
For a connected network of $N$ nodes with system dynamics in (\ref{eq:dt}), we assume that all nodes follow the confidential interaction protocol illustrated in Fig. \ref{fig:1}. If a node Alice is connected to the rest of the network only through an (or a group of colluding) honest-but-curious node(s) Eve, Alice's initial state can be uniquely inferred by Eve in an asymptotic sense.
\end{Theorem 4}

\textit{Proof}: If Alice is directly connected to multiple honest-but-curious nodes that collude with each other, then these nodes can share information with each other to cooperatively estimate Alice's initial state, and hence can be regarded as one node. Therefore, we can only consider the case in which Alice is only connected to one honest-but-curious node Eve, as illustrated in Fig. \ref{fig:privacy}(b). In this case, from the perspective of the honest-but-curious node Eve, the measurement seen at each time step $k$ is $\Delta x_{EA}[k]=a^{(k)}_{A \shortrightarrow E} a^{(k)}_{E \shortrightarrow A} (x_A[k]-x_E[k])$ with $\Delta x_{EA}[k]=-\Delta x_{AE}[k]$. From the system dynamics in (\ref{eq:ex_dt}), Eve knows
\begin{equation}\label{Eve_infer_1}
\begin{aligned}
x_A[0]& = x_A[P] - \varepsilon \sum_{k=0}^{P-1} \Delta x_{AE}[k]\\
& = x_A[P] + \varepsilon \sum_{k=0}^{P-1} \Delta x_{EA}[k]\\
\end{aligned}
\end{equation}
and hence can construct the following observer to estimate the initial value of $x_A[0]$:
\begin{equation}\label{Eve_estimator}
\hat{x}_A[0] = x_E[P]+ \varepsilon \sum_{k=0}^{P-1} \Delta x_{EA}[k]
\end{equation}
As $P$ goes to infinity, average consensus will be achieved asymptotically, i.e., $\lim_{P \shortrightarrow \infty} x_A[P] = \lim_{P \shortrightarrow \infty} x_E[P] = \alpha = \sum_{i=1}^{N} x_i[0] / N$, which leads to
\begin{equation}\label{Eve_infer_2}
\begin{aligned}
\lim_{P \shortrightarrow \infty} \hat{x}_A[0] &= \lim_{P \shortrightarrow \infty} \big(x_E[P] + \varepsilon\sum_{k=0}^{P-1} \Delta x_{EA}[k] \big)\\
& = \lim_{P \shortrightarrow \infty} \big(x_A[P] + \varepsilon\sum_{k=0}^{P-1} \Delta x_{EA}[k] \big)\\
& =x_A[0]
\end{aligned}
\end{equation}
Therefore Eve can uniquely infer Alice's initial state $x_A[0]$ in an asymptotic sense through the estimator (\ref{Eve_estimator}). \hfill{$\blacksquare$}

Theorem \ref{thm:4} implies that the single neighbor configuration should always be avoided, which is also required by other noise-based privacy protocols such as \cite{Mo17} and \cite{Manitara13}.

\subsection{Security Solution}
Due to the additive homomorphic property, the Paillier cryptosystem is vulnerable to active adversaries who may attempt to alter the message being sent through the channel. Although such adversaries cannot find out the exact states of the communicating nodes, they can still inflict significant damage to the system.

\begin{figure}
    \centering
    \includegraphics[width=0.37\textwidth]{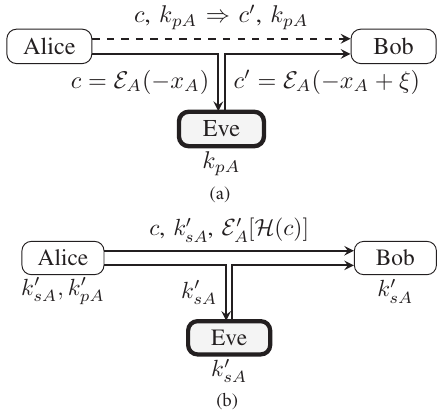}
    \caption{Illustration of attacks from an active attacker (a) and the defense mechanism with a digital signature (b).}\label{fig:security}
\end{figure}

Consider the scenario where the communication from node Alice to Bob is intercepted by an active adversary Eve (cf. Fig. \ref{fig:security}(a)). Since Alice's public key $k_{pA}$ is sent along with $\mathcal{E}(-x_1)$, Eve may use the additive homomorphism to inject an arbitrary noise $\xi$ to the original message $\mathcal{E}(-x_1)$ to sway it to $\mathcal{E}(-x_1+\xi)$. If Bob has no way to tell whether the received message has been modified, Eve may exploit this vulnerability to make the network either converge to a wrong value or not converge at all.

In applications where security is of prime concern, it is imperative to be able to verify the integrity of any received message. Here we propose a solution to provide resilience against active adversaries. Our basic idea is to attach a digital signature to exchanged messages in the confidential interaction protocol, as illustrated in Fig. \ref{fig:security}(b). The signature requires an additional pair of public/private keys ($k'_{pA},k'_{sA}$) and a hash function $\mathcal{H}(\cdot)$, and is represented as ($k_{sA}'$, $\mathcal{E}'_A\left[\mathcal{H}(c)\right]$) where $c$ is the encrypted message. The additional private key $k'_{sA}$ is sent so that Bob can decrypt $\mathcal{E}'_A\left[\mathcal{H}(c)\right]$ and check if the resulting $\mathcal{H}(c)$ matches the received $c$ in terms of the hash operation $\mathcal{H}(\cdot)$. Based on the attached digital signature, the recipient can verify possible modifications during communication.

\begin{Lemma 2}
    When Alice and Bob use the mechanism in Fig. \ref{fig:security}(b) for message transmission and integrity verification, Eve cannot fake Alice's identity by forging a signature that could pass the verification of Bob.
\end{Lemma 2}

\textit{Proof}:
From the defense mechanism illustrated in Fig. \ref{fig:security}(b), it can be seen that without the public key $k'_{pA}$, Eve cannot forge a valid signature (that can be decrypted by Bob), any Eve's attempt to modify $c$ will cause a mismatch between received $c$ and decrypted $\mathcal{H}(c)$ in terms of the hash operation $\mathcal{H}(\cdot)$ (cf. \cite{goldreich2009foundations} for details). Therefore, the defense mechanism can provide resilience against active adversaries. \hfill{$\blacksquare$}

\section{Extensions to Other Consensus}\label{sec:application}

Using the same confidential interaction protocol, we can ensure the privacy of other variants of average consensus. Here we show the applications to three other commonly used consensus problems, i.e., the weighted average consensus, maximum consensus, and minimum consensus.

\subsection{Weighted Average Consensus}

Weighted average consensus seeks convergence of all states to a weighed sum of the initial state, i.e., $\frac{\sum_{j=1}^N w_j x_j[0]}{\sum_{j=1}^N w_j}$ with $w_j>0$ being the weights. According to \cite{Olfati-Saber2007}, weighted average consensus can be achieved by using the following update rule:
\begin{equation} \label{eq:wt_avg_dt}
 x_i[k+1] = x_i[k] + \frac{\varepsilon}{w_i}\sum_{v_j\in N_i} a^{(k)}_{ij} (x_j[k]- x_i[k])
\end{equation}

Note that the average consensus is in fact a special case of the weighted average consensus with all $w_i$ being equal to the same constant.

\begin{Theorem 6}
For a connected network of $N$ nodes, if the coupling weights $a^{(k)}_{ij}$ in (\ref{eq:wt_avg_dt}) are established according to the confidential interaction protocol illustrated in Fig. \ref{fig:1} and the admissible range $[\barbelow{a}, \, \bar{a}]$ for the random selection of $a_{i \shortrightarrow j}^{(k)}$ satisfies $0 < \barbelow{a} < \bar{a} < \sqrt{\frac{\min_{i} w_i}{\varepsilon \Delta}}$, then the system will achieve weighted average consensus with states converging to
    \begin{equation}
    \lim_{k\to \infty} x_i[k]= \frac{\sum_{j=1}^{N} w_jx_j[0]}{\sum_{j=1}^{N} w_j}, \quad i=1,2,\ldots,N.
    \end{equation}
\end{Theorem 6}

\textit{Proof}: The proof can be obtained by following a similar line of reasoning of Theorem 1 in \cite{zhang2011distributed}. \hfill{$\blacksquare$}

\subsection{Maximum Consensus}

The maximum consensus seeks the convergence of all states to the maximal initial value among all states, i.e., $\lim_{k\shortrightarrow \infty}x_i[k]= \max_j x_j[0]$ for $i=1,2,\ldots,N$. Next we show how to use the confidential interaction protocol in Sec. III to guarantee the privacy in maximum consensus. In contrast to the conventional maximum consensus protocol where the current state of a node is directly replaced with the maximum state among its neighbors \cite{maximum_consensus}, we formulate maximum consensus problem as a nonlinear dynamical system to fit into our confidential interaction protocol.

We design the update rule as
\begin{equation} \label{eq:max_con_dt}
 x_i[k+1] = x_i[k] + \max_{v_j\in N_i\cup \{v_i\}} a_{j \shortrightarrow i}^{(k)} (x_j[k]-x_i[k])
\end{equation}

Compared with the average consensus protocol, the difference is that we replaced the summation operation with the ``max'' operator. What is more important is the inclusion of the node itself in the computation, which ensures that the output of the max operator is non-negative.

In this new maximum consensus approach, the weighted state difference can still be calculated using the confidential interaction protocol in Sec. III, but the operation of multiplying a node's own random weight to the weighted difference, i.e., $a_{1 \shortrightarrow 2}$ in (\ref{eq:mul_a1}), can be removed to simplify the computation because there is no need to guarantee a symmetric interaction graph any more. For this reason, here we always set $a_{1 \shortrightarrow 2}$ in (\ref{eq:mul_a1}) to 1, which leads to the coupling weight in (\ref{eq:max_con_dt}) above. Furthermore, to guarantee the convergence of maximum consensus, the admissible range $[\barbelow{a}, \, \bar{a}]$ for the random selection of $a_{j \shortrightarrow i}^{(k)}$ should satisfy $0 < \barbelow{a} < \bar{a} < 1$.

\begin{Theorem 8}\label{thm:8}
For a connected network of $N$ nodes, under the confidential interaction protocol illustrated in Fig. \ref{fig:1} and $a_{j \shortrightarrow i}^{(k)}\in [\barbelow{a}, \, \bar{a}]$ for $k \geq 0$ with $0 < \barbelow{a} < \bar{a} < 1$, the update rule (\ref{eq:max_con_dt}) can achieve maximum consensus, i.e.,
\begin{equation}
\lim_{k\to \infty} {x}_i[k] = \max_{1\leq j\leq N} x_j[0], \quad i=1,2,\ldots,N.
\end{equation}
\end{Theorem 8}

\textit{Proof}: According to (\ref{eq:max_con_dt}), we have
\begin{equation}
\begin{aligned}
 x_i[k+1] - x_i[k] &= \max_{v_j\in N_i\cup \{v_i\}} a_{j \shortrightarrow i}^{(k)} (x_j[k]-x_i[k])\\
 &\geq a_{i \shortrightarrow i}^{(k)} (x_i[k]-x_i[k])= 0
\end{aligned}
\end{equation}
So $x_i[k+1] - x_i[k] \geq 0$ always holds for $i=1,2,\ldots,N$, i.e., the value of each node is non-decreasing.

Let $x_p[k]=\max_i x_i[k]$ at step $k$. According to (\ref{eq:max_con_dt}) and $0 < a_{j \shortrightarrow p}^{(k)} <1$, we have
\begin{equation}
\begin{aligned}
 x_p[k+1]-x_p[k] = \max_{v_j\in N_p\cup \{v_p\}} a_{j \shortrightarrow p}^{(k)} (x_j[k]-x_p[k])\leq 0\\
\end{aligned}
\end{equation}
Since the value of each node is non-decreasing, we have $x_p[k+1]=x_p[k]$, i.e., $x_p[k]$ is time-invariant.

Next we show that all the other nodes will stay less or equal to $x_p[k]$. For every node $i$, we have
\begin{equation} \label{max_con_1}
\begin{aligned}
 x_i[k+1]&= x_i[k] +\max_{v_j\in N_i\cup \{v_i\}} a_{j \shortrightarrow i}^{(k)} (x_j[k]-x_i[k])\\
\end{aligned}
\end{equation}
Since $\max_{v_j\in N_i\cup \{v_i\}} a_{j \shortrightarrow i}^{(k)} (x_j[k]-x_i[k]) \geq 0$ holds, we have
\begin{equation} \label{max_con_2}
\begin{aligned}
& \max_{v_j\in N_i\cup \{v_i\}} a_{j \shortrightarrow i}^{(k)} (x_j[k]-x_i[k]) \leq \max_{v_j\in N_i\cup \{v_i\}} (x_j[k]-x_i[k]) \\
& \leq \max_{v_j\in V} (x_j[k]-x_i[k]) \leq x_p[k]-x_i[k]
\end{aligned}
\end{equation}
where we used the facts $0< a_{j \shortrightarrow i}^{(k)}<1$ and $x_p[k]=\max_i x_i[k]$. Combing (\ref{max_con_2}) with (\ref{max_con_1}) leads to
\begin{equation} \label{max_con_3}
\begin{aligned}
x_i[k+1] & = x_i[k] +\max_{v_j\in N_i\cup \{v_i\}} a_{j \shortrightarrow i}^{(k)} (x_j[k]-x_i[k])\\
& \leq x_i[k] + x_p[k]-x_i[k] = x_p[k]
\end{aligned}
\end{equation}
As a result, the maximum value is invariant with respect to $k$.

Let $\max_i x_i[k] = \max_i x_i[0]=\alpha$. We can now define the error vector as $\boldsymbol{\delta}[k]=\alpha \mathbf{1}-\mathbf{x}[k]$, where $\mathbf{x}[k]=[x_1[k], x_2[k],\cdots, x_N[k]]^T$. Note that ${\delta}_i[k]\geq 0$ holds for all $i$ and $k$. Theorem \ref{thm:8} is thus equivalent to proving
\begin{equation}
 \lim_{k\rightarrow \infty}\boldsymbol{\delta}[k]=\mathbf{0}
\end{equation}

According to (\ref{eq:max_con_dt}), the dynamics of ${\delta}_i[k]$ is given by
\begin{equation}
 {\delta}_i[k+1]-{\delta}_i[k] = -\max_{v_j\in N_i\cup \{v_i\}} a_{j \shortrightarrow i}^{(k)} (x_j[k]-x_i[k])
\end{equation}

Define the Lyapunov function $V(\boldsymbol{\delta}[k]) = \boldsymbol{\delta}[k]^T \boldsymbol{\delta}[k] \geq 0$, where the equality holds only when $\boldsymbol{\delta}[k]=\mathbf{0}$, then we have
\begin{equation}
\begin{aligned}
 V(\boldsymbol{\delta}[k+1])-V(\boldsymbol{\delta}[k])&=\boldsymbol{\delta}[k+1]^T\boldsymbol{\delta}[k+1]-\boldsymbol{\delta}[k]^T\boldsymbol{\delta}[k]
\end{aligned}
\end{equation}
Expanding the right-hand-side (RHS) yields:
\begin{equation}
\begin{aligned}
 & {\text{RHS}} =(\boldsymbol{\delta}[k+1]-\boldsymbol{\delta}[k])^T(\boldsymbol{\delta}[k+1]+\boldsymbol{\delta}[k])\\
 &=-\sum_{v_i\in V}\max_{v_j\in N_i\cup \{v_i\}} a_{j \shortrightarrow i}^{(k)} (x_j[k]-x_i[k]) ({\delta}_i[k+1]+{\delta}_i[k])\\
 &\leq 0
\end{aligned}
\end{equation}
The equality holds when $\max_{v_j\in N_i\cup \{v_i\}} a_{j \shortrightarrow i}^{(k)} (x_j[k]-x_i[k]) = 0$ is true for $i=1,2,\ldots,N$, i.e., when the maximum consensus is achieved. So the Lyapunov function $V(\boldsymbol{\delta}[k])$ will keep decreasing until the maximum consensus is achieved. \hfill{$\blacksquare$}

\subsection{Minimum Consensus}

In contrast to maximum consensus, the minimum consensus problem seeks the convergence of all states to the minimal initial value among all states, i.e., $\lim_{k\shortrightarrow \infty}x_i[k]= \min_j x_j[0]$ for $i=1,2,\ldots,N$. Following the same idea as maximum consensus, we propose a new update rule for achieving minimum consensus to fit into our confidential interaction protocol:
\begin{equation} \label{eq:min_con_dt}
x_i[k+1] = x_i[k] + \min_{v_j\in N_i\cup \{v_i\}} a_{j \shortrightarrow i}^{(k)} (x_j[k]-x_i[k])
\end{equation}

Similar to Theorem \ref{thm:8}, we can obtain the following results
for minimum consensus:

\begin{Theorem 10}\label{thm:10}
For a connected network of $N$ nodes, under the confidential interaction protocol illustrated in Fig. \ref{fig:1} and $a_{j \shortrightarrow i}^{(k)}\in [\barbelow{a}, \, \bar{a}]$ for $k \geq 0$ with $0 < \barbelow{a} < \bar{a} < 1$, the update rule (\ref{eq:min_con_dt}) can achieve minimum consensus, i.e.,
 \begin{equation}
 \lim_{k\to \infty} {x}_i[k] = \min_{1\leq j\leq N} x_j[0], \quad i=1,2,\ldots,N.
 \end{equation}
\end{Theorem 10}

\textit{Proof}: The proof follows the same line of reasoning as Theorem \ref{thm:8}
and is omitted.\hfill{$\blacksquare$}

\section{Implementation Details}\label{sec:impl}
In addition to the constraint imposed on $a_{i \shortrightarrow j}^{(k)}$, there are other technical issues that must be addressed for the implementation of our confidential interaction protocol.

\subsection{Quantization}
Real-world applications typically have $x_i\in \mathbb{R}$ which are represented by floating point numbers in modern computing architectures. On the contrary, encryption algorithms only work on unsigned integers. Define the casting function $f(\cdot, \cdot): \mathbb{R}\times \mathbb{R}\rightarrow \mathcal{M}\subset \mathbb{Z}$ and its inverse $f^{-1}(\cdot, \cdot):\mathcal{M}\times \mathbb{R}\rightarrow \mathbb{R}$ as
\begin{equation}\label{largeQ}
\begin{aligned}
 &f(x, Q) = \left\lceil Qx \right\rfloor_\mathcal{M},\quad
 &f^{-1}(y, Q) = \frac{y}{Q}
\end{aligned}
\end{equation}
where $\left\lceil \cdot \right\rfloor_\mathcal{M}$ maps the input to the nearest integer in $\mathcal{M}$. For the Paillier cryptosystem, this mapping is equivalent to the rounding operation, hence the step size is $1$ which is uniform. Consequently the maximum quantization error is bounded by
\begin{equation}
\begin{aligned}
 \max_{x\in \mathbb{R}}|x-f^{-1}\big( f(x, Q) , Q\big)|= \frac{1}{Q}
\end{aligned}
\end{equation}

Although the quantization error is unavoidable, it can be made arbitrarily small by using a large enough $Q$ in (\ref{largeQ}). In practice we choose a sufficiently large value for $Q$ so that the quantization error is negligible. This is exactly how we convert the state $x_i$ of a node from real value to a fixed length integer and back to a floating point number. The conversion is performed at each iteration of the protocol.

\subsection{Subtraction and Negative Values}
Another issue is how to treat the sign of an integer for encryption. \cite{homomorphic_privacy1} solves this problem by mapping negative values to the end of the group $\mathbb{Z}_n$ where $n=pq$ is given by the public key. We offer an alternative solution by taking advantage of the fact that encryption algorithms blindly treat bit strings as an unsigned integer. In our implementation all integer values are stored in fix-length integers (i.e., \texttt{long int} in C) and negative values are left in two's complement format. Encryption and intermediate computations are carried out as if the underlying data were unsigned. When the final message is decrypted, the overflown bits (bits outside the fixed length) are discarded and the remaining binary number is treated as a signed integer which is later converted back to a real value.

\section{Numerical Simulations and Hardware Experiments}\label{sec:example}

In this section, we first numerically verify effectiveness of the proposed approach and show its advantage over existing state-of-the-art approaches. Then we provide hardware experimental results on a resource-constrained Raspberry-Pi board based micro-controller network to demonstrate the efficiency of the proposed approach.

\subsection{Numerical Implementation of the Proposed Approach}

To illustrate the capability of our protocol, we implemented the consensus protocol in C/C++. We used an open-source C implementation of the Paillier cryptosystem \cite{libpaillier} because it allowed byte-level access. For each exchange between two nodes, the states were converted to 64-bit integers by multiplying $Q=10^5$. The weights $a_{i \shortrightarrow j}^{(k)}$ were also scaled up similarly and represented by 64-bit random integers. The encryption/decryption keys were set to 256-bit long.

\begin{figure}
 \centering
 \includegraphics[width=0.5\textwidth]{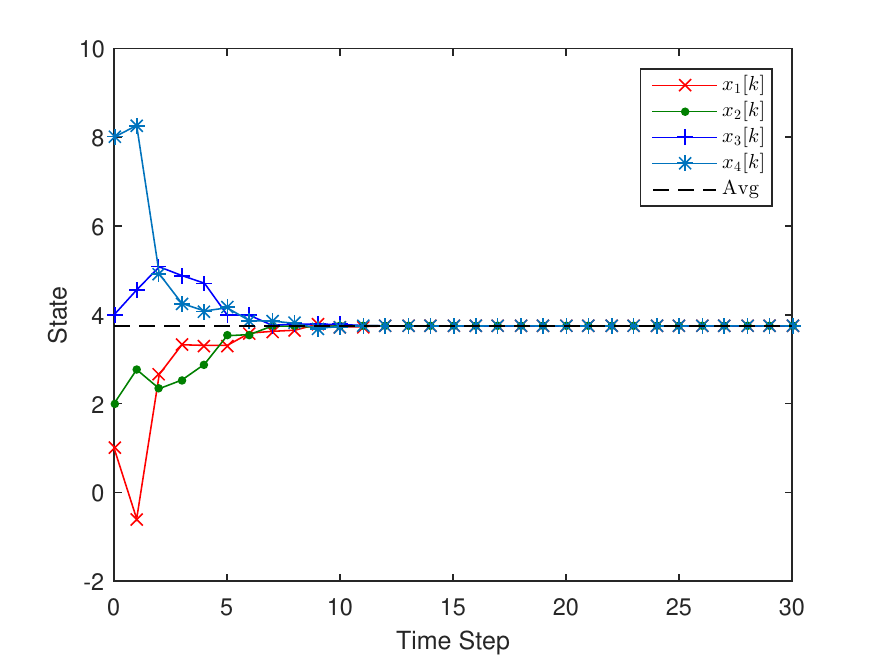}
 \caption{Convergence to the average consensus. The states converge to the average consensus value 3.75 in about 10 steps.}
 \label{fig:eg1_plot}
\end{figure}

\begin{figure}
 \centering
 \includegraphics[width=0.5\textwidth]{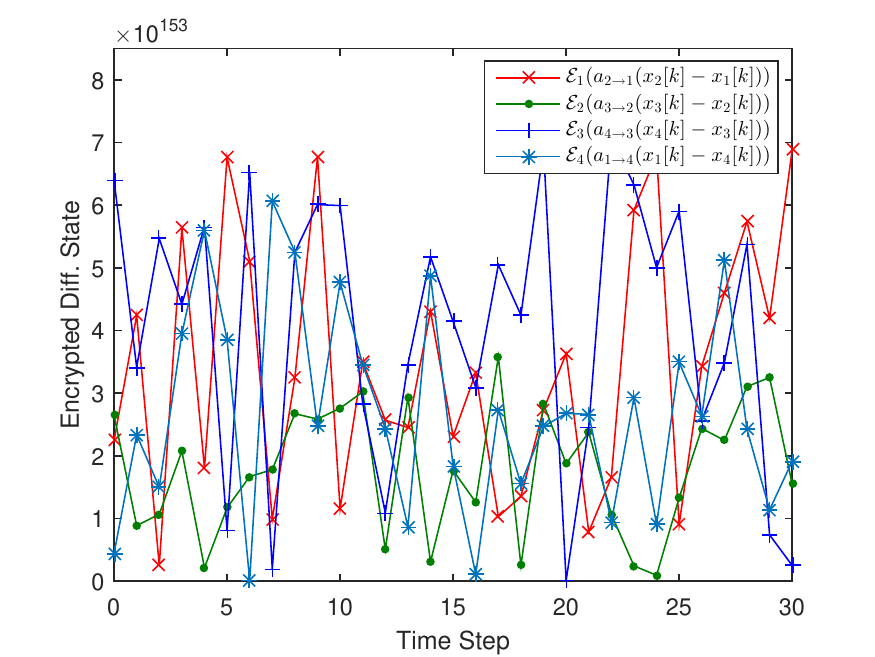}
 \caption{Encrypted weighted differences vs time step. Although the states have converged after 10 steps, the encrypted differences still appear to be random.}
 \label{fig:eg1_edx}
\end{figure}

The first simulation had four nodes connected in a undirected ring, which led to $\Delta=\max_i|N_i|=2$. We set the step-size to $\varepsilon=0.5$. Given $0 < \barbelow{a} < \bar{a} < \frac{1}{\sqrt{\varepsilon \Delta}}=\frac{1}{\sqrt{0.5 *2 }}=1$, we set the admissible range to $[\barbelow{a}, \, \bar{a}]=[0.01, \, 0.99]$. The initial states were set to $\{1, \, 2,\, 4, \, 8\}$ respectively and the average is 3.75. Each node used a static key pair which was initialized once at the beginning. The states' convergence to the average is shown in Fig. \ref{fig:eg1_plot}.

The plot of received encrypted messages which encoded the weighted difference between two nodes is given in Fig. \ref{fig:eg1_edx}. It is worth noting that although the states have converged to the average, the encrypted weighted differences still appear to be random to an unintended observer.

\begin{figure}
 \centering
 \includegraphics[width=0.5\textwidth]{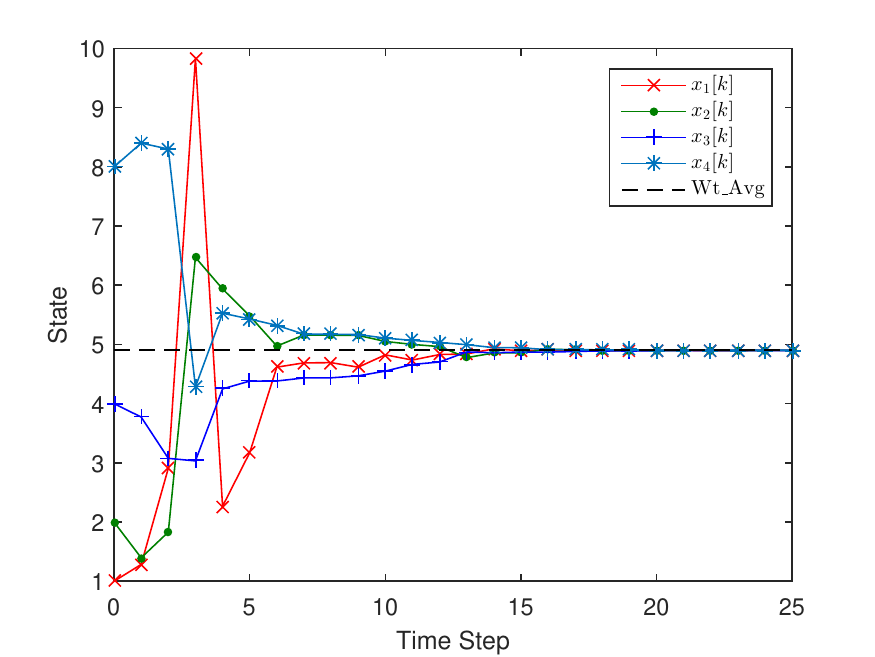}
 \caption{Convergence to the weighted average consensus. The states converge to the weighted average consensus value 4.9 in about 15 steps.}
 \label{fig:eg2_plot}
\end{figure}

The second simulation considered the weighted average consensus. Under the same initial condition as the first example, the nodes also had the associated weights $\{0.1,0.2,0.3,0.4\}$. The step size $\varepsilon$ was set to $0.05$. Under $0 < \barbelow{a} < \bar{a} < \sqrt{\frac{\min_{i} w_i}{\varepsilon \Delta}}= \sqrt{\frac{0.1}{0.05 *2}}=1$, we set the admissible range to $[\barbelow{a}, \, \bar{a}]=[0.01, \, 0.99]$. Fig. \ref{fig:eg2_plot} shows that the states converge to the weighted average value 4.9.

\begin{figure}
 \centering
 \includegraphics[width=0.5\textwidth]{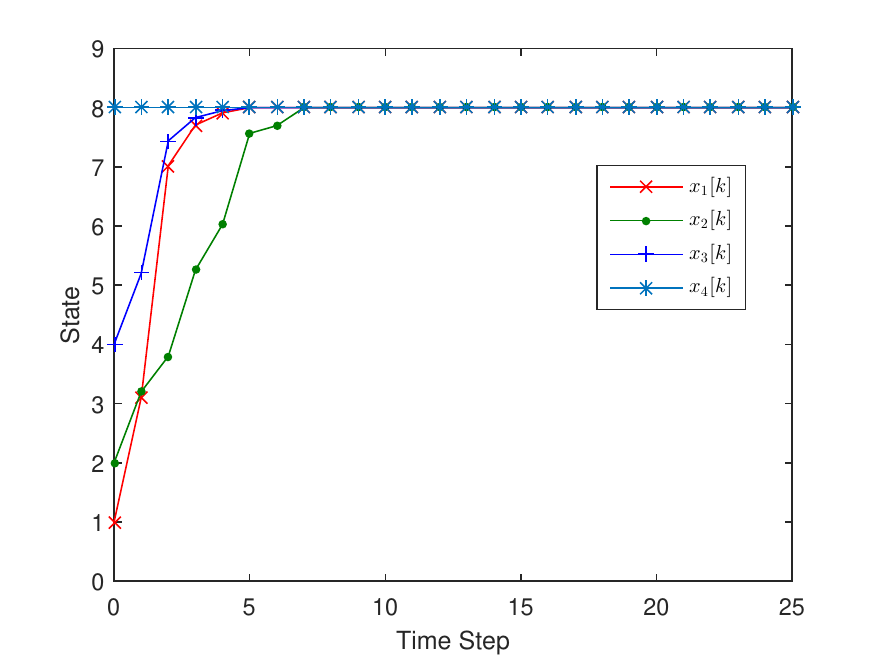}
 \caption{Convergence to the maximum consensus. The states converge to the maximum consensus value 8 in less than 10 steps.}
 \label{fig:eg3_plot}
\end{figure}

The third simulation considered the maximum consensus, using the same graph structure and initial condition. Fig. \ref{fig:eg3_plot} shows that the states converge to the exact maximum in about 7 steps.

The computational overhead caused by the encryption is manageable. Without any hardware-specific optimization, it takes about 7 ms to compute one exchange of state on a desktop computer with a 3.4 GHz CPU and 4.00 GB memory.

\subsection{Comparison with Existing Results}

In this subsection, we show the advantages of the proposed average consensus approach over existing data-obfuscation based average consensus approaches, more specifically, the decaying-noise approach by Mo and Murray \cite{Mo17}, the finite-noise-sequence approach by Manitara and Hadjicostis \cite{Manitara13}, and the differential-privacy based approach by Huang \textit{et al.} \cite{huang2012differentially}. The differential-privacy based approach in \cite{huang2012differentially} injects uncorrelated Laplace noise to exchanged states to ensure the privacy of participating nodes and is subject to a trade-off between privacy and accuracy. To guarantee the accuracy of average consensus, the approaches in \cite{Mo17} and \cite{Manitara13} employ correlated noises, which, however, compromise the achievable privacy. In this subsection, we use numerical simulations to show that all these approaches are vulnerable to an external eavesdropper.

In the numerical simulation, we used the network topology and weight matrix provided in \cite{Mo17}, i.e.,

\begin{equation}
\mathbf{A} = \frac{1}{4}\begin{bmatrix}
2 & 1 & 0 & 0 & 1\\
1 & 2 & 1 & 0 & 0\\
0 & 1 & 2 & 0 & 1\\
0 & 0 & 0 & 3 & 1\\
1 & 0 & 1 & 1 & 1\\
\end{bmatrix}
\end{equation}

We set the initial states of five nodes to $\left\{ 1, \, 2, \, 3, \, 4, \, 5 \right\}$ respectively.

We consider the case where node $1$'s communication with all its neighbors is tapped by an external eavesdropper Eve. Eve also knows the internal dynamics and network interaction topology, i.e., the $\mathbf{A}$ matrix.

\subsubsection{Comparison with the approach in \cite{Mo17}}
We first simulated the approach in \cite{Mo17}. For consistency and ease of reading, we borrow the notations from \cite{Mo17}, where the authors denote the internal state and the transmitted state by $x_i[k]$ and $x_i^+[k]$, respectively, i.e.,
\begin{equation}
x^+_i[k] = x_i[k] + w_i[k]
\end{equation}
where $w_i[k]$ is a random noise specified as follows
\begin{equation}\label{mo_noise}
w_i[k]=
\begin{cases}
v_i[0], & \text{if } k=0\\
\rho^kv_i[k]-\rho^{k-1}v_i[k-1] & \text{otherwise}
\end{cases}
\end{equation}
In (\ref{mo_noise}), $\rho \in (0, \, 1)$ holds and $v_i[k]$ is a normally distributed random variable with zero mean and unit variance. $\{v_i[k]\}_{i=1,\ldots,N, \, k=0,1,\ldots}$ are jointly independent but from (\ref{mo_noise}), it is clear that the injected noise $w_i[k]$ is time-correlated.

Each node updates its internal state by
\begin{equation}
x_i[k+1] = a_{ii}x^+_i[k] + \sum_{v_j\in N_i} a_{ij}x^+_j[k]
\end{equation}

Based on the assumption that Eve has access to the in-bound and out-bound messages of node 1 (i.e., $x_1^+[k]$ and $x_j^+[k]$ for $v_j\in N_1$) as well as the interaction dynamics/topology (i.e., $a_{1j}$), we can build an observer to estimate node 1's initial state as follows:

\begin{itemize}
 \item The initial state of the observer is set to the first out-going transmission from node 1, i.e.
 \begin{equation}
 z[0] = x^+_1[0]
 \end{equation}
 \item The update of the observer is based on accumulating the difference between the transmitted state and the predicted state
 \begin{equation}
 z[k+1]=z[k]+x^+_1[k+1]-\big(a_{11}x^+_1[k] + \sum_{v_j\in N_1} a_{1j}x^+_j[k]\big)
 \label{eq:obs_update}
 \end{equation}
\end{itemize}

\begin{figure}
    \centering
    \includegraphics[width=0.5\textwidth]{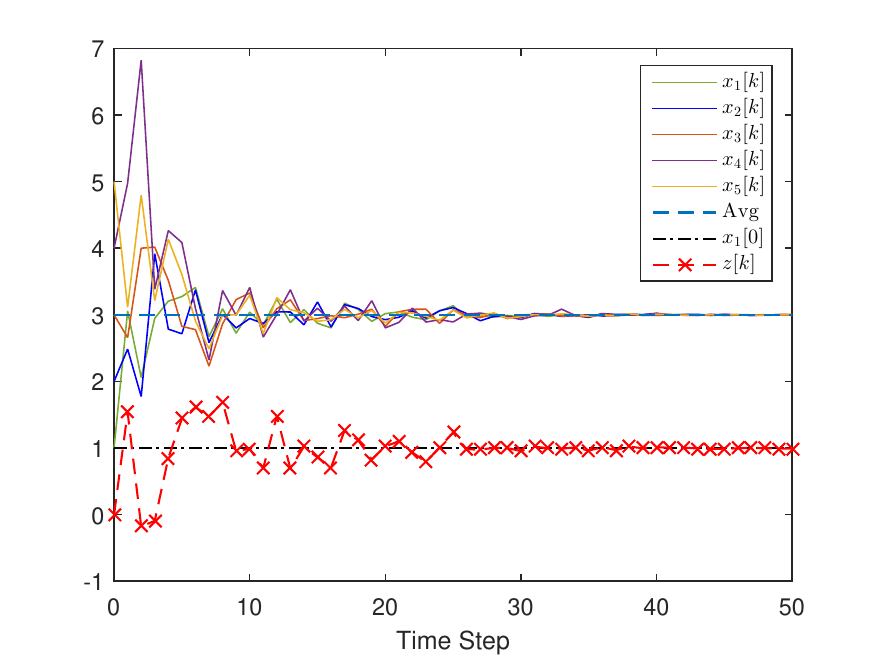}
    \caption{If Eve wiretaps all the messages of node 1 and has access to node dynamics and interaction topology, then the initial state of node 1 can be estimated under the privacy protocol in \cite{Mo17}.}
    \label{fig:mo_obs}
\end{figure}

It is shown in Fig. \ref{fig:mo_obs} that as the network converges and all node states approach the average consensus value, the observer also converges to node 1's true initial value.

In contrast, with our method all the messages are encrypted by design. Consequently, Eve cannot gain any advantage by wiretapping the communication.

\subsubsection{Comparison with the approach in \cite{Manitara13}}
\begin{figure}
 \centering
 \includegraphics[width=0.5\textwidth]{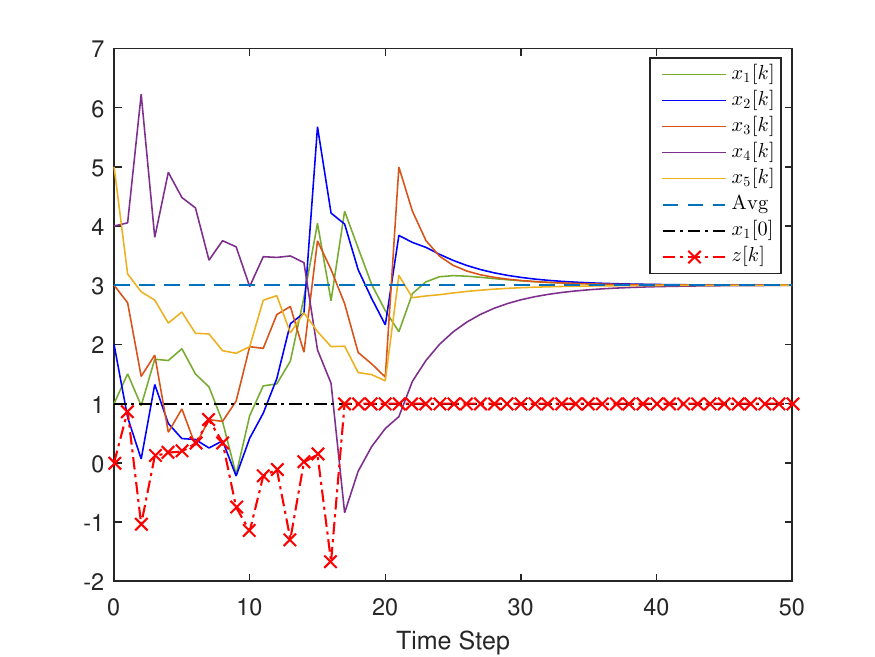}
 \caption{If Eve wiretaps all the messages of node 1 and has access to node dynamics and interaction topology, then the initial state of node 1 can be estimated under the privacy protocol in \cite{Manitara13}.}
 \label{fig:manitara_obs}
\end{figure}

Although using a different noise sequence, the approach in \cite{Manitara13} builds on a similar noise-injection approach to that of \cite{Mo17} (in the sense that later noise will cancel out previous noise). Therefore, we used the same observer design as in the simulation of \cite{Mo17}. Again we obtained that as the network converged, the observer value converged to node $1$'s true initial state, as depicted in Fig. \ref{fig:manitara_obs}.

On the contrary, such external eavesdropper gains no information from our protocol where transmitted messages are encrypted inherently.

\subsubsection{Comparison with the approach in \cite{huang2012differentially}}
We next show the vulnerability of the approach in \cite{huang2012differentially} to an external eavesdropper. We built a naive estimator that guesses the initial state by taking the value of the initial out-bound message of the node. Using the notations from \cite{huang2012differentially}, the estimated state is given by
\begin{equation}
z = x_1[0] = \theta_1[0] + \eta_1[0]
\end{equation}
where $x_1$ is the sent message, $\theta_1[0]$ is the true initial state, and $\eta_1[0]$ is a random noise generated from a Laplace distribution $Lap(cq^0)$.

\begin{figure}
 \centering
 \includegraphics[width=0.5\textwidth]{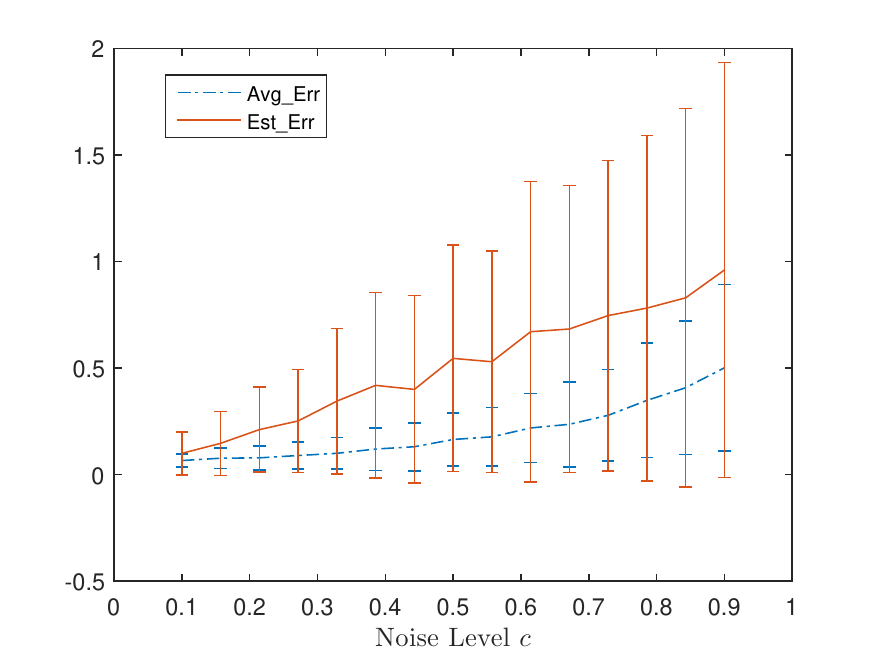}
 \caption{Although the differential-privacy based approach guarantees that the initial state cannot be estimated exactly, the accuracy of the final consensus value and the accuracy of the estimated initial value are correlated.}
 \label{fig:huang_Nfold}
\end{figure}

We then varied the noise parameter $c$ and repeated the simulation for 500 times. The result is summarized in Fig. \ref{fig:huang_Nfold} with Avg\_Err and Est\_Err denoting error in consensus value achieving and initial state estimation respectively. We can conclude that although the differential-privacy based approach can use a large noise level to prevent an external eavesdropper from accurately inferring the initial state, the noise will also lead to a large error in the final consensus value. As a result, when an application calls for higher accuracy of the consensus result, the risk of disclosing one's initial state also becomes higher.

Our method protects the privacy by an encrypted message-exchange mechanism that does not affect the accuracy of the final consensus value. Therefore the user does not have to worry about the trade-off between accuracy and privacy.

\subsection{Hardware Experiments}

\begin{figure}
    \begin{center}
        \includegraphics[width=0.37\textwidth]{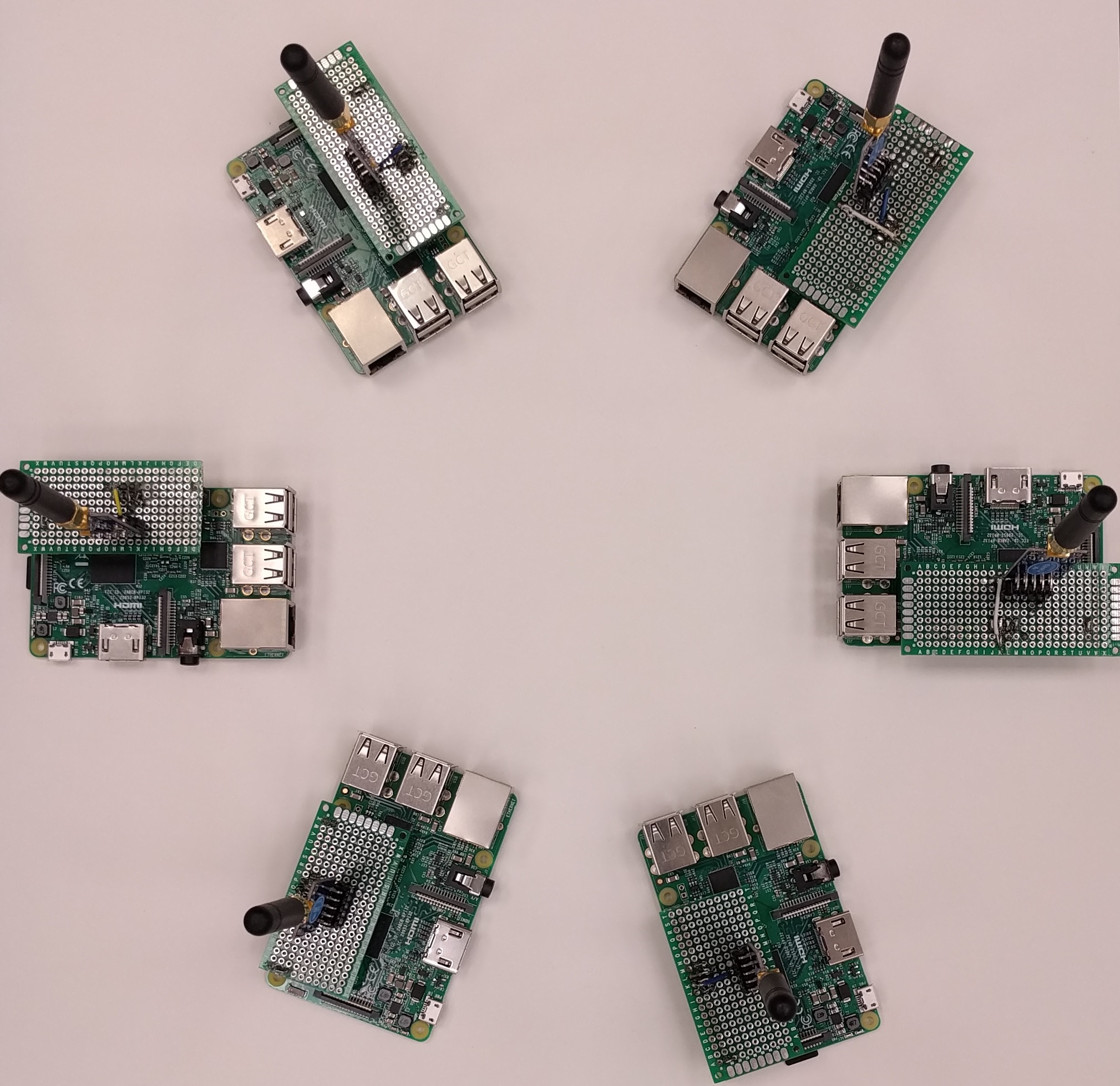}
    \end{center}
    \caption{A network of six Raspberry Pi boards.}
    \label{fig:Raspberry}
\end{figure}

To confirm the efficiency of the secure and privacy-preserving average consensus approach in real-world cyber-physical systems, we also implemented the proposed approach on six Raspberry Pi boards with 64-bit ARMv8 CPU and 1 GB RAM (cf. Fig. \ref{fig:Raspberry}, source code available at \cite{Github:17}).

\begin{figure}
    \begin{center}
        \includegraphics[width=0.5\textwidth]{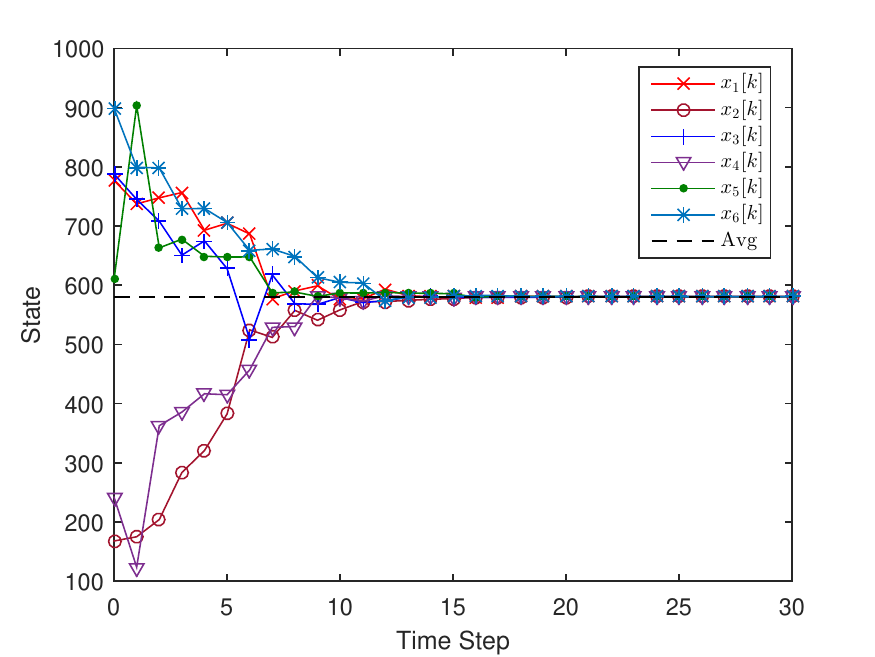}
    \end{center}
    \caption{All nodes converge to the average consensus value in the experimental verification using Raspberry Pi boards. The states have initial values as 777, 168, 788, 242, 610, and 899, respectively and they all converge to the average consensus value 580.67 in about 13 steps.}
    \label{fig:experimental_plot}
\end{figure}

In the implementation, the communication was conducted through Wi-Fi based on the ``sys/socket.h'' C library. Paillier encryption/decryption was realized using the ``libpaillier-0.8'' library from \cite{Paillier_lib}. To obtain $\Delta x_{ij}$ in a pairwise interaction, a node employs a request message to initialize the interaction and the other node replies with a response message. In a multi-node network, for a node to be able to simultaneously receive requests and responses from multiple neighbors, parallelism needed to be introduced. The ``pthread'' C library was used to generate multiple parallel threads to handle incoming requests and responses. Each time a node receives a request/response, it generates a new thread to handle it and immediately listens for more requests. Because in the implementation, it is impossible to start all nodes simultaneously, a counter is introduced on each node and its value is embedded in each request/response packet to help nodes make sure that they are on the same pace. For 64 byte encryption key, the size of the actual packet is 144 bytes, which includes all necessary headers and stuffing bytes. For each interaction, the average processing latency was 7.8 ms, which is acceptable for most real-time control systems. The implementation result is given in Fig. \ref{fig:experimental_plot}, which shows that perfect consensus can be achieved.

\section{Discussions and Conclusions}\label{sec:conclusion}
In this paper we proposed a decentralized secure and privacy-preserving protocol for the network consensus problem which can guarantee the security and privacy of a node as long as it has at least one legitimate neighbor. In contrast to existing data-obfuscation based average consensus approaches which are vulnerable to an external eavesdropper, we encode randomness into the system dynamics with the help of an additive homomorphic cryptosystem which allows a deterministic convergence to the average (subject to a quantization error). The protocol also provides resilience to passive attackers and allows easy incorporation of active attacker defending mechanisms. In addition to average consensus, our protocol can be easily extended to enable security and privacy in weighted average consensus and maximum/minimum consensus under new update rules. Both numerical simulations and hardware experiments are provided to demonstrate the effectiveness and efficiency of the approach.

The proposed approach offers several advantages over existing privacy-preserving average consensus protocols using noise-injection based data obfuscation, such as the uncorrelated-noise (differential-privacy) based approach in \cite{huang2012differentially} and the correlated-noise based approaches in \cite{Mo17,Manitara13}. First, compared with the differential-privacy based approach, our approach can guarantee convergence to the consensus value (subject to a quantization error) in a deterministic manner. Secondly, our approach encrypts the exchanged messages and is resilient to external eavesdroppers wiretapping the communication channel, whereas an external eavesdropper can easily defeat the approaches in \cite{Mo17} and \cite{Manitara13} by using a naive observer, as confirmed by the numerical simulation results in Sec. VIII-B. Furthermore, the proposed approach is naturally extendable to time-varying networks, whereas existing noise-injection based data obfuscation approaches all assume time-invariant parameters, which could be troublesome if the topology of the network or the number of nodes is not constant.

On the other hand, although the computational overhead introduced by the encryption algorithm (7 ms in numerical simulations and 7.8 ms in hardware experiments) is acceptable in most real-time control systems, it is indeed higher compared to the unencrypted alternatives. We argue that the benefits of using encryption to preserve privacy and security outweigh the computational burden which is easily manageable on modern micro-controllers.

\section*{Acknowledgement}
The authors would like to thank Christoforos Hadjicostis and Yilin Mo for their comments on an initial draft of this article.

\bibliographystyle{IEEEtran}
\bibliography{reference1}

\end{document}